\documentclass[11pt]{article}
\usepackage{eurosym}
\usepackage{amsfonts}
\usepackage{amsmath}
\usepackage{geometry}
\usepackage{amssymb}
\usepackage{graphicx}
\usepackage{color}
\usepackage{tikz}
\usepackage[author={Oana}]{pdfcomment}
\usepackage[colorinlistoftodos]{todonotes}
\usepackage{titling}
\usepackage{lipsum}
\usepackage{mathrsfs} 

\newtheorem{theorem}{Theorem}

\newtheorem{definition}[theorem]{Definition}

\newtheorem{lemma}[theorem]{Lemma}

\newtheorem{proposition}[theorem]{Proposition}
\newtheorem{remark}[theorem]{Remark}

\begin{document}

\date{}
\title{\textbf{Well-posedness for a stochastic 2D Euler equation with
transport noise}}
\author{ Oana Lang \thanks{Corresponding author, email address: o.lang15@imperial.ac.uk .}\qquad Dan Crisan \\
{\small Department of Mathematics, Imperial College London, London, SW7 2AZ, UK}}
\maketitle

\begin{abstract}
We prove the existence of a unique global strong solution for a stochastic two-dimensional
Euler vorticity equation for incompressible flows with noise of transport type. In particular, we show that the initial smoothness of the solution is preserved.  The arguments are based on approximating the solution of the Euler
equation with a family of viscous solutions which is proved to be relatively
compact using a tightness criterion by Kurtz \cite{EthierKurtz}.  \\
\\


\end{abstract}

\vspace{7mm}

\section{Introduction}

Consider the two dimensional Euler equation modelling an incompressible flow perturbed by transport type
stochasticity 
\begin{equation*}
d{\omega }_{t}+u_{t}\cdot \nabla {\omega }_{t}dt+\displaystyle%
\sum_{i=1}^{\infty }\mathsterling_i {\omega }_{t}\circ dW_{t}^{i}=0
\end{equation*}%
with initial condition $\omega _{0}$, where $\xi _{i}$ are time-independent divergence-free vector fields and $(W^{i})_{i \in \mathbb{N}}$ is a sequence of independent Brownian motions.
Classically, $u_{t}$ stands for the velocity of an incompressible fluid and $%
\omega _{t}=curl\ u_{t}$ is the corresponding fluid vorticity. The stochastic part
considered here follows the Stochastic Advection by Lie Transport (SALT) theory (see \cite%
{Holm2015}, \cite{CFH}, \cite{FFnoisetransport}, \cite{Wei1}) and corresponds to a stochastic integral of Stratonovich type. \\

The Euler equation is used to model the motion of an incompressible inviscid fluid. A representative aspect in this context is the study of the fluid vortex dynamics modelled by the vorticity equation. There is a vast literature on well-posedness in the deterministic setting, see e.g. \cite{MajdaBertozzi}, \cite{Judovich63}, \cite{EbinMarsden}, \cite{Temam}, \cite{Kato}, \cite{KatoLai}, \cite{Bardos}, \cite{BourguignonBrezis}, and references therein.

 The introduction of stochasticity into ideal fluid dynamics has received special attention over the past two decades. On one hand, comprehensive physical models can be obtained when the stochastic term accounts for physical uncertainties (\cite{Holm2015}, \cite{Wei1}, \cite{Wei2}, \cite{CFH}), whilst, in some cases, the regularity properties of the deterministic solution can be improved when the right type of stochasticity is added (\cite{FFnoisetransport}, \cite{FlandoliFedrizzi}, \cite{VicolHoltz}, \cite{Delarue}). Global existence of smooth solutions for the stochastic Euler equation with multiplicative noise in both 2D and 3D has been obtained in \cite{VicolHoltz}. In \cite{BessaihFerrario}, a weak solution of the Euler equation with additive noise is constructed as an inviscid limit of the stochastic damped 2D Navier-Stokes equations. A martingale solution constructed also as a limit of Navier-Stokes equations but with cylindrical noise can be found in \cite{BrzezniakPeszat}. Existence and uniqueness results with different variations in terms of stochastic forcing and approximations can be found in \cite{MikuleviciusValiukevicius1}, \cite{CapinskiCutland}, \cite{MikuleviciusValiukevicius2}, \cite{Bessaihmartingale99}, \cite{RocknerLiu} and references therein. An overview of results on this topic is provided in \cite{Bessaih2013}. 
 
The analysis of nonlinear stochastic partial differential equations with noise of
transport type has recently expanded substantially, see e.g., \cite{Holm2015}, \cite{CFH}, \cite{Wei1}, \cite{Wei2}, 
\cite{Bendall1}, \cite{Aythami1}, \cite{Aythami2}, \cite{GayBalmaz}.
Existence of a solution  for the two-dimensional
stochastic Euler equation with noise of
transport type has been considered in \cite{Brzezniak}. While in \cite{Brzezniak} the authors  prove the existence and pathwise uniqueness of a distributional  solution in $L^{\infty}(%
\mathbb{T}^2)$, in this paper we are concerned with the existence of a strong solution and give conditions under which the solution enjoys  smoothness properties\footnote{ In other 
words, we identify conditions under which the (strong) solution of the two-dimensional
stochastic Euler equation with noise of
transport type belongs to the Sobolev space $\mathcal{W}^{k,2}$ with $k$ arbitrarily high.}. In \cite{FlandoliLuo}, a random point vortices system is
used to construct a so-called $\rho$ - white noise solution. Local
well-posedness and a Beale-Kato-Majda blow-up criterion for the
three-dimensional case in the space $\mathcal{W}^{2,2}(\mathbb{T}^3)$ has
been obtained in \cite{CFH}. Full well-posedness for a point vortex dynamics
system based on this equation has been proven in \cite{Flandolivortex}. The linear case has been considered in \cite{FlandoliFedrizzi} and in \cite{FlandoliGubinelliPriola}. 

In the sequel, $\mathbb{T}^2$ is the two-dimensional torus, $k\ge 2$ is a fixed positive integer and $\mathcal{W}^{k,2}$ is the usual Sobolev space (see Section \ref{preliminaries}). 
The main result of this paper is the following: \\[2mm]
\textit{Theorem: Under certain conditions on the vector fields $(\xi_i)_i$ the two-dimensional stochastic Euler vorticity equation 
\begin{equation}\label{maineqn}
d{\omega }_{t}+u_{t}\cdot \nabla {\omega }_{t}dt+\sum_{i=1}^{\infty }(\xi
_{i}\cdot \nabla {\omega }_{t})\circ dW_{t}^{i}=0,\ \ \  \omega
_{0}\in \mathcal{W}^{k,2}(\mathbb{T}^{2}), \ \ \ 
\end{equation}%
admits a unique global (in time) solution which belongs to the space $\mathcal{W}^{k,2}(\mathbb{T}^{2})$. Moreover, $\omega_t$ is a continuous function of the initial condition.}

\begin{remark}\label{itoform}
As stated above, the stochastic terms in \eqref{maineqn} are  stochastic integrals of Stratonovich type. We interpret equation \eqref{maineqn} in its corresponding It\^{o} form, that is
\begin{equation}\label{maineqnito}
d\omega_t + u_t \cdot \nabla \omega_tdt + \displaystyle\sum_{i=1}^{\infty} \xi_i \cdot \nabla\omega_tdW_t^i = \frac{1}{2}\displaystyle\sum_{i=1}^{\infty} \xi_i \cdot \nabla\big(\xi_i \cdot \nabla \omega_t\big)dt.
\end{equation}
\end{remark}

The assumptions on the vector fields $(\xi_i)_i$ are described in Section \ref{preliminaries}. In short, they are assumed to be sufficiently smooth, their corresponding norms to decay sufficiently fast as $i$ increases, so that the infinite sums in \eqref{maineqn}, respectively, in \eqref{maineqnito} make sense in the right spaces (see condition \eqref{xiassumpt} below). Importantly, we do not require the additional assumption\footnote{Here $c$ is a non-negative constant, $I_2$ is the identity matrix, and $\xi(x)^{\star}$ is the transpose of $\xi(x)$.} 
\begin{equation}\label{adas}
\sum_{i=1}^{\infty} \xi_i(x)\xi_i^{\star}(x) = cI_2 \ \   
\end{equation}
used in \cite{Brzezniak}. As a result, in the It\^{o} version \eqref{maineqnito} of the SPDE, the term $\frac{1}{2}\displaystyle\sum_{i=1}^{\infty} \xi_i \cdot \nabla(\xi_i \cdot \nabla \omega_t)$ does not reduce to  $c\Delta\omega_t$. This would simplify the analysis as, in this case, the Laplacian commutes with higher order derivatives. Morever, it commutes with the operation of convolution with the Biot-Savart kernel, an essential ingredient used   in \cite{Brzezniak}. The general term $\frac{1}{2}\displaystyle\sum_{i=1}^{\infty} \xi_i \cdot \nabla(\xi_i \cdot \nabla \omega_t)$ makes the analysis harder. We succeed in controlling it by considering it in tandem with the term $\displaystyle\sum_{i=1}^{\infty}\displaystyle\int_{\mathbb{T}^2} (\xi_i \cdot \nabla\omega_t)^2dx$ coming from the quadratic variation of the stochastic integrals (see Lemma \ref{apriori} i.) appearing in the evolution equation for the process $t\mapsto\|\omega_t\|^2$. A similar technical difficulty appears when trying to control the high-order derivatives of the vorticity. 
Nonetheless, this is achieved through a set of inequalities (see Lemma \ref{apriori}) that have first been introduced in the literature by Krylov and Rozovskii (\cite{KrylovRozovskii}, \cite{GerencserGyongyKrylov}) and recently used by Crisan, Flandoli and Holm (\cite{CFH}).  Again, we emphasize that these rather surprising inequalities hold true without imposing assumptions on the driving vectors $(\xi_i)_i$ other than on  their smoothness and summability. This finding is particularly important when using this model for the purpose of uncertainty quantification and data assimilation: for example,  in  \cite{Wei1}, \cite{Wei2}, \cite{Wei3}, the driving vectors $(\xi_i)_i$  are \emph{estimated} from data and not a priori chosen. The methodology used in these papers does not naturally lead to driving vector fields that satisfy assumption \eqref{adas} so removing it is essential for our research programme.

We emphasize that the appearance of the second order differential operator  $\omega\mapsto\frac{1}{2}\sum_{i=1}^{\infty} \xi_i \cdot \nabla(\xi_i \cdot \nabla \omega)$ in the It\^o version of the Euler equation does not give the equation a parabolic character, even if one assumes the restriction \eqref{adas} with $c$ chosen strictly positive. Equation \eqref{maineqn} is truly a transport type equation and one cannot expect the initial condition to be smoothed out. The best scenario is to prove that the initial level of smoothness of the solution is preserved. This is indeed the main finding of our research. Moreover, we show that our result can be extended to cover also $L^{\infty}$-solutions in the Yudovich sense.      

The paper is organised as follows: In Section \ref{preliminaries} we introduce the main assumptions, key notations, and some preliminary results. In Section \ref{mainresults} we present our main results: in subsection \ref{uniquenessEuler} we show that the solution is almost surely pathwise unique, while in subsection \ref{existenceEuler} we prove existence of a strong solution (in the sense of Definition \ref{solutions}). In Section \ref{existenceEulertruncated} we proceed with an extensive analysis of a truncated form of the Euler equation: uniqueness (Section \ref{uniquenessEulert1}) and existence - based on a new approximating sequence introduced in Section \ref{existenceEulert}. At the end of this section we show continuity with respect to initial conditions for the original equation. In Section \ref{approx.sequence} we show existence, uniqueness, and continuity for the approximating sequence of solutions constructed in Section \ref{existenceEulert}.  In Section \ref{sectiontightness} we show that the family of approximating solutions is relatively compact. In Section \ref{yudovichsection} we present an extension of the main results to the Yudovich setting. The paper is concluded with an Appendix that incorporates a number of proofs of the technical lemmas and statements of some classical results. 

\section{Preliminaries}\label{preliminaries}

We summarise the notation used throughout the manuscript. Let $(\Omega, \mathcal{F}, (\mathcal{F}_t)_t, \mathbb{P})$ be a filtered probability space, with the sequence $\left(W^{i}\right)_{i\in\mathbb{N}}$ of
independent Brownian motions defined on it. Let $X$ a generic Banach space. Throughout the paper $C$ is a generic notation for constants whose values can change from line to line.

\begin{itemize}
\item We denote by $\mathbb{T}^{d}=\mathbb{R}^{d}/\mathbb{Z}^{d}$ the d-dimensional
torus. In our case $d=2$. 
\item Let $\alpha = (\alpha_1, \alpha_2, \ldots, \alpha_d) \in \mathbb{N}^d, \ d>0$ be a multi-index of length $|\alpha| = \displaystyle\sum_{j=1}^{d}\alpha_j$. Then $\partial^{\alpha} = \partial^{(\alpha_1, \alpha_2, \ldots, \alpha_d)} = \partial_1^{\alpha_1}\partial_2^{\alpha_2} \ldots \partial_d^{\alpha_d}$ denotes the differential operator of order $|\alpha|$, with $\partial^{(0,0, \ldots, 0)}f = f$ for any function $f$ defined on $\mathbb{T}^d$, and $\partial_i^{\alpha_i} = \frac{\partial^{\alpha_i}}{\partial x_i^{\alpha_i}}$, $x \in \mathbb{T}^d$. In our case $d=2$ and $|\alpha| \leq k$. 
\item $L^{p} (\mathbb{T}^{2}; X)$ \footnote{Here and later whenever the space $X$ coincides with the Euclidean space $\mathbb{R}$ or $\mathbb{R}^2$, it is omitted from the notation: For example $L^{p} (\mathbb{T}^{2}; X)$ becomes $L^{p} (\mathbb{T}^{2})$, etc.} is the class of all measurable $p$ - integrable functions $f$ defined on the two-dimensional torus, with values in $X$ ($p$ is a positive real number). The space is endowed with its canonical norm $\|f\|_p = \bigg(\displaystyle\int_{\mathbb{T}^2} \|f\|_X^p dx\bigg)^{1/p}.$ Conventionally, for $p = \infty$ we denote by $L^{\infty}$ the space of essentially bounded measurable functions. 
\item For $a,b\in L^{2}\left(  \mathbb{T}^{2}\right)$, we denote by $\langle \cdot ,\cdot \rangle $ the
scalar product
\begin{equation*}
\langle a,b\rangle :=\displaystyle\int_{\mathbb{T}^{2}}a(x)\cdot b(x)dx.
\end{equation*}%
\item $\mathcal{W}^{m,p}(\mathbb{T}^{2})$ is the Sobolev space of functions $f \in L^p(\mathbb{T}^2)$ such that $D^{\alpha}f \in L^p(\mathbb{T}^2)$ for $0 \leq |\alpha| \leq m$, where $D^{\alpha}f$ is the distributional derivative of $f$.
The canonical norm of this space is $\|f\|_{m,p} = \bigg(\displaystyle\sum_{0 \leq |\alpha| \leq m} \|D^{\alpha}f\|_p^p\bigg)^{1/p}$, with $m$ a positive integer and $1\leq p < \infty$. 
A detailed presentation of Sobolev spaces can be found in \cite{Adams}. 
\item $C^m(\mathbb{T}^2; X)$ is the (vector) space of all $X$-valued functions $f$ which are continuous on $\mathbb{T}^2$ with continuous partial derivatives $D^{\alpha}f$ of orders $|\alpha| \leq m$, for $m \geq 0$. $C^{\infty}(\mathbb{T}^2; X)$ is regarded as the intersection of all spaces $C^{m}(\mathbb{T}^2; X)$. Note that on the torus all continuous functions are bounded. 
\item $C([0, \infty); X)$ is the space of continuous functions from $[0, \infty)$ to $X$ equipped with 
the uniform convergence norm over compact subintervals of $[0, \infty)$.
\item 
$L^p(0, T; X)$ is the space of measurable functions from $[0, T]$ to $X$ such that the norm 
\[
\|f\|_{L^p(0, T; X)}=\bigg(\displaystyle\int_0^T \|f(t)\|_X^pdt\bigg)^{1/p}
\] 
is finite. 
 
\item $D([0, \infty); X)$ is the space of c\`{a}dl\`{a}g functions, that is functions $f:[0, \infty) \rightarrow X$ which are right-continuous and have limits to the left, endowed with the Skorokhod topology. This topology is a natural choice in this case because its underlying metric transforms $D([0, \infty); X)$ into a complete separable metric space. For further details see \cite{EthierKurtz} Chapter 3, Section 5, pp. 117-118. 
\item Given $a: \mathbb{T}^{2}\rightarrow\mathbb {R}^2$, we define the differential
operator $\mathsterling_a$ by 
$\mathsterling_{a}b:=a\cdot \nabla b$
 for any map $b: \mathbb{T}^{2}\rightarrow\mathbb {R}$ such that the scalar product
 $a\cdot \nabla b$ makes sense. In line with this, we use the notation
\begin{equation*} 
\mathsterling_{i}\omega _{t}:=\mathsterling_{\xi _{i}}\omega _{t}:=\xi
_{i}\cdot \nabla \omega _{t}\ \ \ \hbox{and}\ \ \ \mathsterling%
_{i}^{2}\omega _{t}:=\mathsterling_{\xi _{i}}^{2}\omega _{t}:=\xi _{i}\cdot
\nabla (\xi _{i}\cdot \nabla \omega _{t}).
\end{equation*}%
Denote the dual of $\mathsterling_{i}$ by $\mathsterling_{i}^{\star }$ that is $%
\langle \mathsterling _{i}a,b\rangle =\langle a,\mathsterling_{i}^{\star }b\rangle$. 
\item For any vector $u \in \mathbb{R}^2$ we denote the gradient of $u$ by $\nabla u = (\partial_1u, \partial_2 u)$ and the corresponding orthogonal by $\nabla^{\perp}u = (\partial_2 u, -\partial_1 u )$.
\end{itemize}

\begin{remark}\label{dual}
If $div \ \xi _{i}=\nabla \cdot \xi _{i}=0$, then the dual of the operator $\mathsterling_{i} $ is $-\mathsterling_{i}$.
\end{remark}
\noindent\textbf{Assumptions on the vector fields $(\xi_i)_i$}. 
The vector fields $\xi_i:\mathbb{T}^2 \rightarrow \mathbb{R}%
^2 $ are chosen to be time-independent and divergence-free
and they need to be specified from the underlying physics. We assume that 
\begin{equation} \label{xiassumpt}
\displaystyle\sum_{i=1}^{\infty} \|\xi_i\|_{k+1,\infty}^2 < \infty.
\end{equation}
Condition \eqref{xiassumpt} ensures that for any $f \in \mathcal{W}^{2,2}(\mathbb{T}^2)$
\begin{subequations}
\begin{equation}\label{xiconda}
   \displaystyle\sum_{i=1}^{\infty} \|\mathsterling_{i}f\|_{2}^2 \leq C\|f\|_{1,2}^2 
\end{equation}
\begin{equation}\label{xicondb}
    \displaystyle\sum_{i=1}^{\infty} \|\mathsterling_i^2f\|_2^2 \leq C\|f\|_{2,2}^2.
\end{equation}
\end{subequations}
Provided $\omega\in L^2(0, T; W^{2,2}(\mathbb{T}^2, \mathbb{R}))$, condition \eqref{xiconda} ensures that the infinite sum of stochastic integrals
\begin{equation} \label{infsum1}
\sum_{i=1}^{\infty} \int_0^t \mathsterling_{i} \omega_sdW_s^i
\end{equation}
is well defined and belongs to $L^2(0, T; L^2(\mathbb{T}^2, \mathbb{R}))$.
Similarly, condition \eqref{xicondb} ensures that the process $s\rightarrow\mathsterling_i^2\omega_s$
is well-defined and belongs to $L^2(0, T; L^2(\mathbb{T}^2, \mathbb{R}))$ provided the solution of the stochastic partial differential equation \eqref{maineqn}  belongs to a suitably chosen space (see Definition \ref{solutions} below). In particular, the It\^o correction in \eqref{maineqnito} is well defined.
The conditions above are needed also for proving a number of required a priori estimates (see Lemma \ref{apriori} in Appendix).
\begin{definition}\label{solutions}$\left.\right.$\\[-5mm]
\begin{enumerate}
\item[a.] A strong solution of the
stochastic partial differential equation \eqref{maineqn} is an $(\mathcal{F}_{t})_t$
- adapted process $\omega :\Omega \times \mathbb{T}^{2}\rightarrow \mathbb{R}$ with trajectories in the space $C([0,\infty);\mathcal{W}^{k,2}(\mathbb{T}^{2}))$, such that the identity\footnote{Here and everywhere else $u$ is implicitly defined as the velocity field whose vorticity is $\omega$, in other words $\omega =curl\ u=\partial_2 u^1-\partial_1 u^2$. See further details in Remarks \ref{whatisu} and \ref{biotsavartkernel}.}
\begin{equation}\label{strongdefinition}
\omega _{t}=\omega _{0}-\int_{0}^{t}\mathsterling_{u_s}  \omega
_{s}ds-\sum_{i=1}^{\infty }\int_{0}^{t}\mathsterling_{i} \omega
_{s}dW_{s}^{i}+\frac{1}{2}\sum_{i=1}^{\infty }\int_{0}^{t}\mathsterling_{i} ^{2}\omega _{s}ds
\end{equation}%
with $\omega _{|_{t=0}}=\omega _{0}$, holds $\mathbb{P}$ - almost surely.
\footnote{Equation 
(\ref{strongdefinition})  is interpreted as an identity between elements in   $L^{2}(\mathbb{T}^{2};\mathbb{R})$. The same applies to the identity (\ref{martingaledefinition}) }
\item[b.] A weak/distributional solution of
equation \eqref{maineqn} is an $(\mathcal{F}_t)_t$ - adapted process $\omega : \Omega
\times \mathbb{T}^2 \rightarrow \mathbb{R}$ with trajectories in the set $%
C([0,\infty); L^2(\mathbb{T}^2,\mathbb{R})) 
$, which satisfies the
equation \eqref{maineqn} in the weak topology of $L^2(\mathbb{T}^2,\mathbb{R})$, i.e. 
\begin{equation}\label{weakdefinition} 
\langle\omega_t, \varphi\rangle = \langle\omega_0,\varphi\rangle + %
\displaystyle\int_0^t \langle \omega_s, \mathsterling_{u_s}^{\star}\varphi
\rangle ds + \displaystyle\sum_{i=1}^{\infty} \displaystyle\int_0^t
\langle\omega_s, \mathsterling_{i}^{\star}\varphi\rangle dW_s^i + \frac{1}{2}%
\displaystyle\sum_{i=1}^{\infty}\displaystyle\int_0^t\langle\omega_s, %
\mathsterling_{i}^{\star}\mathsterling_{i}^{\star}\varphi\rangle ds
\end{equation}
holds $\mathbb{P}$ - almost surely for all $\varphi\in C^{\infty}(\mathbb{T}%
^2, \mathbb{R})$.
\item[c.] A martingale solution of equation \eqref{maineqn} is a triple $(\check\omega, (\check W^i)_i), (\check\Omega, \mathcal{\check F}, \check{\mathbb{P}}), (\mathcal{\check F}_t)_t$ such that $(\check\Omega, \mathcal{\check F}, \check{\mathbb{P}})$ is a probability space, $(\mathcal{\check F}_t)_t$ is a filtration defined on this space,   
$\check\omega$ is a continuous $(\mathcal{\check F}_t)_t$-adapted  real valued process $\check\omega : \Omega
\times \mathbb{T}^2 \rightarrow \mathbb{R}$ with trajectories in the set $%
C([0, \infty); \mathcal{W}^{k,2}(\mathbb{T}^2))$, $(\check W^i)_i$ are independent $(\mathcal{\check F}_t)_t$-adapted Brownian motions and the identity 
\begin{equation*}\label{martingaledefinition}
\check{\omega} _{t}=\check{\omega} _{0}-\int_{0}^{t}\mathsterling_{\check{u}_s}  \check{\omega}
_{s}ds-\sum_{i=1}^{\infty }\int_{0}^{t}\mathsterling_{i} \check{\omega}
_{s}d\check{W}_{s}^{i}+\frac{1}{2}\sum_{i=1}^{\infty }\int_{0}^{t}\mathsterling_{i} ^{2}\check{\omega}_{s}ds
\end{equation*}%
with $\check{\omega} _{|_{t=0}}=\check{\omega}_{0}$, holds $\check{\mathbb{P}}$-almost surely.\footnote{We use the "check" notation in the description of the various components of a weak probabilistic solution, to emphasize that the existence of a weak solution does not guarantee that, for a \underline{given} set of Brownian motions $(W^i)_i$ defined on a (possibly different) probability space $(\Omega, \mathcal{F}, {\mathbb{P}})$ a solution of \eqref{maineqn} will exist. Clearly the existence of a strong solution implies the existence of a martingale solution.
}

\item[d.]A classical solution of equation \eqref{maineqn}
is an $(\mathcal{F}_t)_t$ - adapted process $\omega : \Omega \times \mathbb{T%
}^2 \rightarrow \mathbb{R}$ with trajectories of class $C([0, \infty); C^2(%
\mathbb{T}^2;\mathbb{R}))$.
\end{enumerate}
\end{definition}
\begin{remark}\label{whatisu} The velocity field $u$ is not uniquely identified through the equation $\omega=curl\  u$. Indeed any two velocity fields that differ by a constant will lead to the same vorticity map $\omega$. Instead we identify $u$ through the "explicit" formula $u=\nabla^{\perp}\Delta^{-1}\omega$, see details in Remark \ref{biotsavartkernel} in the Appendix. In particular, since $u$ and $\omega$ are defined in terms of partial derivatives of other
functions, on the torus, they must have zero average:  
\[
\int_{\mathbb{T}^{2}}u^{1}\left(  x\right)  dx=\int_{\mathbb{T}^{2}}u^{2}\left(  x\right)  dx=\int_{\mathbb{T}^{2}}\omega\left(  x\right)  dx=0.
 \]   
This is due to the fact that $\omega = curl \ u = \partial_2u^1 - \partial_1 u^2$ therefore $\displaystyle\int_{\mathbb{T}^2} \omega dx = \displaystyle\int_{\mathbb{T}^2}\partial_2u^1 - \partial_1 u^2 dx = 0$ as we have periodic boundary conditions. Similarly, $u=\nabla^{\perp} \psi $ where $\psi$ is the streamfunction (see the Appendix) and therefore $\displaystyle\int_{\mathbb{T}^2} u dx = 0$.
Note that if $\omega_0$ has zero average, then $\omega_t$ will have zero average, as it is immediate that all the terms appearing in the Euler equation (either \eqref{maineqn} or \eqref{maineqnito}) have zero average.          
\end{remark} 
\begin{remark}
Note that $\omega_t \in \mathcal{W}^{k,2}(\mathbb{T}^2)$ implies $u_t \in \mathcal{W}^{k+1,2}(\mathbb{T}^2)$ (see the Appendix). By standard Sobolev embedding theorems $\mathcal{W}^{k+1,2}(\mathbb{T}^2) \hookrightarrow\mathcal{W}^{k,2}(\mathbb{T}^2) \hookrightarrow L^\infty(\mathbb{T}^2)$ for $k\geq 2$, hence the terms $\mathsterling_{u_t}\omega_t = u_t \cdot \nabla \omega_t\in L^2(\mathbb{T}^2, \mathbb{R})$ in 
\eqref{strongdefinition}, and $\langle \omega_s, \mathsterling_{u_s}^{\star}\varphi\rangle$ in
\eqref{weakdefinition} are well defined. 
However, to ensure that a \emph{classical} solution (Definition 3.d) exists, we require $u\in C^2({\mathbb T}) $ and therefore we need $k\ge 4$.
\end{remark}
\begin{remark} \label{weak} Naturally, if $\omega_t$ is a strong solution in
the sense of Definition \ref{solutions}, then it is also a\\ weak/distributional solution.
In this sense, our result enhances the solution properties presented in \cite%
{Brzezniak} at the expense of stronger assumptions on the initial condition of the stochastic partial differential equation, but without the need to impose the additional constraint \eqref{adas}.  Note also that if $\omega_t$ is a weak/distributional solution with paths in $L^2([0,T]; \mathcal{W}^{k,2}(\mathbb{T}^2))$ then, by integration by parts, the equation has a strong solution. 
\end{remark}

\section{Main results}\label{mainresults}

 We restate the existence and uniqueness result, this time with complete details: 
\begin{theorem} \label{mainresult1}
If $\omega_0 \in \mathcal{W}^{k,2}(\mathbb{T}^2)$, then the two-dimensional stochastic Euler vorticity equation \eqref{maineqn}
\begin{equation*}
d{\omega}_t+ \mathsterling_{u_t}{\omega}_t dt + \displaystyle %
\sum_{i=1}^{\infty}\mathsterling_{i} ^{} \omega_t\circ dW_t^i=0  \label{1}
\end{equation*}
admits a  unique global $(\mathcal F_t)_t$-adapted strong solution
$\omega=\{\omega_t,t\in [0,\infty)\}$ with values in the space $C\left( [  0,\infty);\mathcal{W}^{k,2}(\mathbb{T}^2)  \right)  $. In particular,
if  $k\geq 4$ the solution is classical. 
\end{theorem}
The proof of Theorem \ref{mainresult1} is contained in Sections \ref{uniquenessEuler} and \ref{existenceEuler}. We state next a result that shows the continuity with respect to initial conditions:    
\begin{theorem}\label{mainresult2}
Let $\omega$, $\tilde\omega$ be two strong solutions of equation \eqref{maineqn}. Define $A$ as the process $A_t:=\displaystyle\int_0^t\|\omega_s\|_{k,2}ds$, for any $ t\ge 0$. Then there exists a positive constant $C$ independent of the two solutions, such that 
\begin{equation}\label{conteqn}
\mathbb E [e^{-CA_t}||\omega_t-\tilde\omega_t||_{k,2}^2]\le 
||\omega_0-\tilde\omega_0||_{k,2}^2.
\end{equation}
\end{theorem}
The proof of Theorem \ref{mainresult2} is incorporated in Section \ref{contEuler}.

 \begin{remark}\label{apptomodelling}
 Insofar as Theorem \ref{mainresult1} and Theorem \ref{mainresult2} are less general than, for example, the corresponding results in \cite{Brzezniak}, since the initial condition of equation \eqref{maineqn} is assumed to be in $W^{k,2}(\mathbb{T}^2)$. The relaxation to initial conditions in $L^{\infty}(\mathbb{T}^2)$ is done in Section \ref{yudovichsection} where the well-posedness is achieved without the additional constraint \eqref{adas}. However, the importance of these results lies in their application to modelling and to the numerical analysis of equation \eqref{maineqn}. In particular, Theorem \ref{mainresult1} states that the initial smoothness of the solution \underline{is carried over for all times}. From a modelling perspective, this is quite important: should the vorticity of a fluid be modelled by equation \eqref{maineqn}, it is essential to have it uniquely  defined everywhere. This is not the case if vorticity is only known to be in $L^\infty(\mathbb{T}^2)$. As an immediate consequence of a Sobolev embedding theorem, this is achieved, for example, if $\omega\in W^{2,2}(\mathbb{T}^2)$. Separately, if equation \eqref{maineqn} is used to model the evolution of the state in a data assimilation problem (see, e.g. \cite{Wei3}
 for further details) and the observable is, say, the fluid velocity at a given set of points, then that observable should be well defined \underline{at those chosen points}. Again just by assuming $\omega\in L^\infty(\mathbb{T})$ does not ensure this. Finally, Theorem \ref{mainresult1} and Theorem \ref{mainresult2} are also important when one is interested in the numerical approximation of equation \eqref{maineqn}. For example, when a finite element numerical approximation is used, the smoothness of the solution influences the choice of the basis and governs the rate of convergence of the numerical approximation (higher rates of convergence require higher smootheness), see \cite{BS}. Moreover, Theorem \ref{mainresult2} is useful to transfer the convergence from local to global of numerical approximations (in the required Sobolev norm).
 \end{remark} 
 \begin{remark}\label{gerencserresult}
At the expense of additional technical arguments,  Theorems \ref{mainresult1} and \ref{mainresult2} can be extended to cover solutions in $W^{m,p}(\mathbb{T}^2)$. More precisely, the existence of solutions in $W^{m,p}(\mathbb{T}^2)$ of the sequence of linearised truncated equations in Section \ref{existenceEulert} can be obtained by an extension (albeit non-trivial)  of the results in \cite{GerencserGyongyKrylov} from deterministic to random coefficients. In addition, the results in Lemmas \ref{correctedlemma} and \ref{apriori} would need to be extended to cover the necessary a priori bounds in $W^{m,p}(\mathbb{T}^2)$  \footnote{These two extensions of Theorems \ref{mainresult1} and \ref{mainresult2} have been suggested to us by an anonymous referee. We thank the referee for this.}.
 \end{remark}
 
\subsection{Pathwise uniqueness of the solution of the Euler equation   \label{uniquenessEuler}}
The uniqueness of the solution of equation \eqref{maineqn} is an immediate corollary of inequality \eqref{conteqn} with $k=0$. However, the proof of \eqref{conteqn} requires the existence of an approximating sequence which is constructed as part of the existence results. We present below a direct proof which avoids this, given the fact that pathwise uniqueness is required for  the proof of existence of a strong (probabilistic) solution. 

Suppose that equation \eqref{maineqn} admits two global $(\mathcal{F}_t)_t$-adapted solutions $\omega_1$ and $%
\omega_2 $ with values in the space 
$C\left( [  0,\infty);\mathcal{W}^{k,2}(\mathbb{T}^2)  \right)  $ and let $\bar{\omega}:=\omega^1-\omega^2$. Consider the
corresponding velocities $u^1$ and $u^2$ such that $curl \ u^1=\omega^1$, $%
curl \ u^2=\omega^2$ and $\bar{u}:=u^1-u^2$. Since both $\omega^1$ and $%
\omega^2$ satisfy \eqref{maineqnito}, their difference satisfies 
\begin{equation*}
d{\bar{\omega}}_t+(\mathsterling_{\bar{u}_t}\omega_t^1 + \mathsterling_{u_t^2} {\bar{\omega}_t}) dt + \displaystyle \sum_{i=1}^{\infty}\mathsterling_{i}{%
\bar{\omega}_t} dW_t^i-\frac{1}{2}\displaystyle\sum_{i=1}^{\infty}\mathsterling_{i}^2 \bar{\omega}_tdt=0.
\end{equation*}
By an application of the It\^{o} formula one obtains 
\begin{equation*}
\begin{aligned}
d\|\bar{\omega}_t\|_2^2&=-2\displaystyle\sum_{i=1}^{\infty}\langle
\bar{\omega}_t, \mathsterling_{i}\bar{\omega}_t \rangle dW_t^i-2\langle
\bar{\omega}_t, \mathsterling_{\bar{u}_t}\omega_t^1 + \mathsterling_{u_t^2} {\bar{\omega}_t}\rangle dt \\ &+ \displaystyle\sum_{i=1}^{\infty}
(\big\langle \bar{\omega}_t, \mathsterling_{i}^2
\bar{\omega}_t \big\rangle  + \langle\mathsterling_{i} \bar{\omega}_t, \mathsterling_{i} \bar{\omega}_t \rangle)
dt. \end{aligned}
\end{equation*}
Note that the first and the last terms in the above identity are null (see Lemma \ref{apriori})\footnote{The application of the Lemma requires that the two solutions $\omega^1$ and $\omega^2 $ belong to $\mathcal{W}^{k,2}(\mathbb{T}^2)$ with $k\ge 2$. To deduce \eqref{conteqn} we need a similar control (albeit not an identity) for higher order derivatives. This is done by using the approximating sequence constructed in Section \ref{existenceEulertruncated} and then taking the limit. This is the reason why we cannot prove directly \eqref{conteqn}.} and that
\begin{equation*}
 |\langle
\bar{\omega}_t, \mathsterling_{\bar{u}_t}\omega_t^1 \rangle| \leq
\|\bar{\omega}_t\|_2\|\bar{u}_t\|_{4}\|\nabla\omega_t^1\|_{4} \leq
C\|\bar{\omega}_t\|_2^2\|\omega_t^1\|_{k,2}.
\end{equation*}
This is true since by the Sobolev embedding theorem (see \cite{Adams} Theorem 4.12 case A) 
one has $\|\nabla\omega_{t}^1\|_{4} \leq C \|\omega_t^1\|_{k,2}$ and using also the Biot-Savart law one has
$\|\bar{u}_t\|_{4} \leq C\|\bar{u}_t\|_{1,2} \leq C \|\bar{\omega}_t\|_{2}$. 
Finally, observe that 
$\langle
\bar{\omega}_t,  \mathsterling_{u_t^2} {\bar{\omega}_t}\rangle =-\frac{1}{2} \displaystyle\int_{\mathbb{T}^2} (\nabla \cdot u_t^2)
(\bar{\omega}_t)^2 dx =0 $
since $div \ u_t^2=0$. It follows that 
\begin{equation}\label{lastinequality}
\begin{aligned} d\|\bar{\omega}_t\|_2^2&=-2\langle \bar{\omega}_t, \mathsterling_{\bar{u}_t}\omega_t^1 \rangle dt \leq
C\|\omega_t^1\|_{k,2}\|\bar{\omega}_t\|_{2}^2dt. \end{aligned}
\end{equation}
Since we only have a priori bounds for the expected value of the process $t\rightarrow \|\omega_t^1\|_{k,2}$ and not for its pathwise values, the uniqueness cannot be deduced through a Gronwall type argument. Instead, we proceed as follows: let $A$ be the process defined as 
$A_t:=\displaystyle\int_0^tC\|\omega_s^1\|_{k,2}ds$, for any $ t\ge 0$.
This is an increasing process that stays finite $\mathbb{P}$-almost surely for all $t\ge 0$ as the paths of $\omega^1$ are in $C\left( [  0,\infty);\mathcal{W}^{k,2}(\mathbb{T}^2)  \right)  $. By product rule, 
\begin{equation*}
d\big(e^{-A_t}\|\bar{\omega}_t\|_{2}^2\big)=e^{-A_t}(d\|\bar{\omega}_t\|_{2}^2-C\|\bar{\omega}_t\|_{2}^2\|%
\omega_t^1\|_{k,2}dt) \le 0.
\end{equation*}
This leads to 
\begin{equation*}
\begin{aligned} e^{-A_t}\|\bar{\omega}_t\|_{2}^2 &=0.
\end{aligned}
\end{equation*}
We conclude that 
$e^{-A_t}\|\bar{\omega}_t\|_{2}^2=0$, and since $e^{-A_t}$ cannot be null due to the (pathwise) finiteness of $A_t$ 
we deduce that $\|\bar{\omega}_t\|_{2}^2=0$ almost surely, which gives the claim.

The above argument uses the fact that the terms  $\big\langle \bar{\omega}_t, \mathsterling_{i}^2
\bar{\omega}_t \big\rangle $ and $\langle\mathsterling_{i} \bar{\omega}_t, \mathsterling_{i} \bar{\omega}_t \rangle$ are well defined. In other words, even though we only wish to control the $L^2$-norm of the vorticity, we have to resort to higher order derivatives. This is permitted as we assumed that $\omega\in\mathcal{W}^{k,2}(\mathbb{T}^2)$ for $k\ge 2$. By applying a similar argument, to control the $\mathcal{W}^{k,2}(\mathbb{T}^2)$-norm of the vorticity we would need to control terms of the form $ \big\langle \partial^k\bar{\omega}_t, \partial^k\mathsterling_{i}^2
\bar{\omega}_t\big\rangle $ and $\langle\partial^k\mathsterling_{i} \bar{\omega}_t, \partial^k\mathsterling_{i} \bar{\omega}_t \rangle$. This is no longer allowed because we do not have enough smoothness in the system. To overcome this difficulty we will make use of a smooth approximating sequence for the vorticity equation, see Section \ref{contEuler}.

\begin{remark}
The above uniqueness result is somewhat stronger than the uniqueness deduced from inequality \eqref{conteqn}. It shows that a solution of \eqref{maineqn} will be unique in the larger space    
$L^{2}(\mathbb{T}^2)$ rather than in the space $\mathcal{W}^{k,2}(\mathbb{T}^2)$. Nevertheless, inequality \eqref{conteqn} shows the continuity of the solution with respect to initial conditions.  
\end{remark}

\begin{remark}
We note that, in contrast to the deterministic version of the Euler equation, the minimal $k$ that ensures the existence of a strong solution is $k=2$. This is because of the occurrence of the term  $\mathsterling_{i} ^{2}\omega$ in the It\^o version of the equation \eqref{maineqn}. Moreover, if we insist on the Stratonovich representation of equation \eqref{maineqn}, then we need to use the evolution equation of $\mathsterling_{i}\omega$ to deduce the covariation between $\mathsterling_{i}\omega$  and $W^i$ required for the rigorous definition of the Stratonovich integral. This, in turn, requires $k\ge 3$ as the term $\mathsterling_{i} ^{3}\omega$ appears in this evolution equation. 

Nonetheless, the methodology from this paper can be used to cover initial conditions $\omega_0 \in \mathcal{W}^{k,2}(\mathbb{T}^2)$ with $k<2$. In this case we have to content ourselves with weak/distributional solutions. Whilst this is not the subject of this paper, such a solution can be shown to exist as long as the product $\omega\mathcal{L}_u\varphi$ makes sense in a suitably chosen sense. We need to interpret the nonlinear term in a weak form as a generalised function and replace it with $\displaystyle\int_0^t \langle  \omega_s, u_s \cdot \nabla\varphi\rangle ds$. Then the same methodology can be applied as long as we can control $u\omega$ in a suitably chosen norm. In section \ref{yudovichsection} we do this for the case  $\omega_0 \in L^\infty(\mathbb{T}^2)$(this is the so-called Yudovich setting). 
\end{remark}

{\subsection{Existence of the solution of the Euler equation}\label{existenceEuler}}

The existence of the solution of equation \eqref{maineqn} is proved by first showing that a truncated version of it has a solution, and then removing the truncation. In particular we will truncate the non-linear term in \eqref{maineqn} by using a smooth function $f_R$ equal to $1$ on $[0,R]$, equal to $0$ on $[R+1, \infty)$, and decreasing on $[R, R+1]$, for arbitrary $R>0$. We then have the following: \\
\begin{proposition} \label{truncatedEulert}
If $\omega_0 \in \mathcal{W}^{k,2}(\mathbb{T}^2)$, then the following equation 
\begin{equation}\label{truncatedEulere}
d{\omega}_t^{R}+ f_R(\|\omega_t^{R}\|_{k-1,2})\mathsterling_{u_t^R}{\omega}_t^R dt + \displaystyle %
\sum_{i=1}^{\infty}\mathsterling_{i} ^{} \omega_t^R\circ dW_t^i=0  
\end{equation}
admits a  unique global $(\mathcal F_t)_t$-adapted solution
$\omega^R=\{\omega^R_t,t\in [0,\infty)\}$ with values in the space \\
$C\left( [  0,\infty);\mathcal{W}^{k,2}(\mathbb{T}^2)  \right)  $. In particular,
if  $k\geq 4$, the solution is classical. 
\end{proposition}
\begin{remark}
The truncation in terms of the norm $\|\omega_t^{R}\|_{k-1,2}$ and not $\|\omega_t^{R}\|_{k,2}$ is not incidental as it suffices to control the norm $\|u_t^{R}\|_{k,2}$  (see Proposition \ref{biotsavartlaw}).    
\end{remark}
We prove Proposition \ref{truncatedEulert} in Section \ref{existenceEulertruncated}.
For now let us proceed with the proof of global existence for the solution of the Euler equation \eqref{maineqn}.
\begin{proposition}\label{localtoglobalprop}
The solution of the stochastic 2D Euler equation \eqref{maineqn} is global. 
\end{proposition}
\noindent\textbf{\textit{Proof}} Define $\tau_R(\omega):=\inf_{t\ge 0}\{\|\omega_t^{R}\|_{k-1,2}\geq R\}$. Observe that on $[0,\tau_R]$, $f_R(\|\omega_t^{R}\|_{k-1,2})=1$, and therefore, on $[0,\tau_R]$ the solution of the truncated equation \eqref{truncatedEulere} is, in fact a solution of \eqref{maineqn} with all required properties. It therefore makes sense to define the process  $\omega=\{\omega_t,t\in [0,\infty)\}$
$\omega_t=\omega_t^R$ for $t\in [0,\tau_R]$. This definition is consistent as, following the uniqueness property of the solution of the truncated equation (see Section \ref{uniquenessEulert1}),  
      $\omega_t^{R} =\omega_t^{R'}$ for $t\in [0,\tau_{\min (R,R')}]$. The process  $\omega$ defined this way is a solution of the Euler equation \eqref{maineqn} on the interval   $[0,\displaystyle\sup_R\tau_R)$. To obtain a global solution we need to prove that   $\displaystyle\sup_{R>0}\tau_R = \infty$ $\mathbb{P}$ - almost surely. Let $\mathscr{A}:=\{\omega \in \Omega| \displaystyle\sup_{R>0}\tau_R(\omega) < \infty\}$. Then 
    
\[
\mathscr{A} = \displaystyle\bigcup_{N>0} \displaystyle\{\sup_R\tau_R(\omega) \le N\}=\displaystyle\bigcup_N \bigcap_R \displaystyle\{|\tau_R(\omega)| \le N\}\]
 and 
        \begin{equation*}
         \mathbb{P}\left(| \tau_R(\omega) |\le N\right) = 
        \mathbb{P} \left(\displaystyle\sup_{t\in [0,N]}\|\omega_t^{R}\|_{k-1,2} > R\right).
        \end{equation*}
In order to finish the proof of global existence we use Lemma \ref{correctedlemma} and the fact that
\begin{equation*}
\mathbb{P}\left( \|\omega_t^R\|_{k-1,2} > R\right) \leq \frac{\mathbb{E}\left[ \ln\left(\|\omega_t^R\|_{k-1,2}^2 + e \right)\right]}{R^2 + e} \leq \frac{\mathcal{C}(\omega_0,T)}{R^2 + e} \xrightarrow[R \rightarrow \infty]{}0.
\end{equation*} 
It follows that
\[
\mathbb{P}\left(\bigcap_R \displaystyle\{|\tau_R(\omega)| \le N\}\right)=\lim_{R\rightarrow \infty }\mathbb{P}\left(| \tau_R(\omega)| \le N\right)=0.
\] 
and therefore $\mathbb{P}(\mathscr{A})=0$. This concludes the global existence for the solution of the equation \eqref{maineqn}.

\section{Analysis of the truncated equation}\label{existenceEulertruncated}

\subsection{Uniqueness of solution for the truncated equation} \label{uniquenessEulert1}
\
We use a similar strategy as the one used to prove the uniqueness of the solution of the (un-truncated) Euler equation \eqref{maineqn}. Suppose that equation \eqref{truncatedEulere} admits two global $(\mathcal F_t)_t$-adapted solutions $\omega^{1,R}$ and $%
\omega^{2,R} $ with values in the space 
$C\left( [  0,\infty);\mathcal{W}^{k,2}(\mathbb{T}^2)  \right)$. We prove that $\omega^{1,R}$ and $\omega^{2,R}$ must coincide. In the following,  we will formally drop the dependence on $R$ of the two solutions. As above, let $\bar{\omega}:=\omega^1-\omega^2$ and consider the
corresponding velocities $u^1$ and $u^2$ such that $curl \ u^1=\omega^1$, $%
curl \ u^2=\omega^2$ and $\bar{u}:=u^1-u^2$. Since both $\omega^1$ and $%
\omega^2$ satisfy \eqref{truncatedEulere}, their difference satisfies 
\begin{equation*}
d{\bar{\omega}}_t+\big(\big(K_R(\omega_t^1)\mathsterling_{u^1_t}-K_R(\omega_t^2)\mathsterling_{u^2_t}\big)\omega_t^1 + K_R(\omega_t^2)\mathsterling_{u_t^2} {\bar{\omega}_t}\big) dt + \displaystyle \sum_{i=1}^{\infty}\mathsterling_{i}{%
\bar{\omega}_t} dW_t^i-\frac{1}{2}\displaystyle\sum_{i=1}^{\infty}\mathsterling_{i}^2 \bar{\omega}_tdt=0.
\end{equation*}
where $K_R(\omega)=f_R(\|\omega\|_{k-1,2})$. By an application of the It\^{o} formula and after eliminating the null terms (see Lemma \ref{apriori}, Remark \ref{truncatedsimilar}, and use the fact that $u_t^2$ is divergence-free), one obtains  
\begin{equation*}
d\|\bar{\omega}_t\|_2^2+2\displaystyle\sum_{i=1}^{\infty}\langle
\bar{\omega}_t, \mathsterling_{i}\bar{\omega}_t \rangle dW_t^i=-2\langle
\bar{\omega}_t,(K_R(\omega_t^1)\mathsterling_{u^1_t}-K_R(\omega_t^2)\mathsterling_{u^2_t})\omega_t^1  \rangle dt.  
 \end{equation*}
One can show that (see \cite{CFH} for a proof) there exists a constant $C=C(R)$ such that  
\[
\|K_R(\omega_t^1)u^1_t-K_R(\omega_t^2){u^2_t}\|_4\le C\|\bar \omega_t\|_{k-1,2}
\]
and to finally deduce that  
\begin{equation*}
 |\langle
\bar{\omega},(K_R(\omega_t^1)\mathsterling_{u^1_t}-K_R(\omega_t^2)\mathsterling_{u^2_t})\omega_t^1  \rangle| 
\leq
C\|\bar{\omega}_t\|_2\|\bar \omega_t\|_{k-1,2}\|\omega_t^1\|_{k,2}. 
\end{equation*}
It follows that (note that the stochastic term is null)
\begin{equation*}
\begin{aligned} d\|\bar{\omega}_t\|_2^2&
\leq C\|\omega_t^1\|_{k,2}\|\bar{\omega}_t\|_{k,2}^2dt. \end{aligned}
\end{equation*}
Similar arguments are used to control $\|\partial^\alpha \bar{\omega}_t\|_2^2$  where $\alpha$ is a multi-index with $|\alpha|\le k-1$ and to deduce that there exists a constant $C$ such that 
\[
\begin{aligned} d\|\partial^\alpha\bar{\omega}_t\|_2^2&+2\displaystyle\sum_{i=1}^{\infty}\langle\partial^\alpha
\bar{\omega}_t, \partial^\alpha\mathsterling_{i}\bar{\omega}_t \rangle dW_t^i\leq
C\|\omega_t^1\|_{k,2}\|\bar{\omega}_t\|_{k,2}^2dt, \end{aligned}
\] 
where we use the control (see Lemma \ref{apriori})
\[
\big\langle \partial^\alpha\bar{\omega}_t, \partial^\alpha\mathsterling_{i}^2
\bar{\omega}_t \big\rangle  + \langle\partial^\alpha\mathsterling_{i} \bar{\omega}_t, \partial^\alpha\mathsterling_{i} \bar{\omega}_t \rangle\le C\|\bar\omega\|_{k,2}^2.
\]
Some care is required for the case when $|\alpha|= k-1$ as $\partial^\alpha\mathsterling_{i}^2
\bar{\omega}_t$ is no longer well-defined. In this case, by using the weak form of the equation \eqref{truncatedEulere} to rewrite $\big\langle \partial^\alpha\bar{\omega}_t, \partial^\alpha\mathsterling_{i}^2\bar{\omega}_t \rangle$ as $-\big\langle \partial^{\alpha_1} \partial^\alpha\bar{\omega}_t, \partial^{\alpha_2}\mathsterling_{i}^2\bar{\omega}_t \rangle$ we can proceed as above by using that 
\[
-\big\langle \partial^{\alpha_1} \partial^\alpha\bar{\omega}_t, \partial^{\alpha_2}\mathsterling_{i}^2\bar{\omega}_t \rangle  + \langle\partial^\alpha\mathsterling_{i} \bar{\omega}_t, \partial^\alpha\mathsterling_{i} \bar{\omega}_t \rangle\le C\|\bar{\omega}_t\|_{k,2}^2.
\]
The above is true for functions in $\mathcal{W}^{k+1,2}(\mathbb{T}^2)$ and, by the continuity of both sides in the above inequality, it is also true for functions in the larger space $\mathcal{W}^{k,2}(\mathbb{T}^2)$, since $\mathcal{W}^{k+1,2}(\mathbb{T}^2)$ is dense in $\mathcal{W}^{k,2}(\mathbb{T}^2)$. The proof is concluded in an identical manner as that for the uniqueness of the Euler equation (see Section \eqref{uniquenessEuler})

\subsection{Existence of solution for the truncated equation} \label{existenceEulert}

The strategy of proving that the truncated equation \eqref{truncatedEulere} has a solution is to construct an approximating sequence of processes that will converge in distribution to a solution of \eqref{truncatedEulere}. This justifies the existence of a weak solution. Together with the pathwise uniqueness of the solution of this equation, we then deduce that strong uniqueness holds.     

Recall that  $\omega_0 \in \mathcal{W}^{k,2}(\mathbb{T}^2)$. Let $(\omega_0^{n})_n \in C^{\infty}(\mathbb{T}^2)$ be a sequence such that $\omega_0^{n} \xrightarrow{n\rightarrow\infty}\omega_0$ in $\mathcal{W}^{k,2}(\mathbb{T}^2)$. 
For any $t\geq 0$ we construct the sequence $(\omega_t^{\nu_n,R,n})_{n\geq 0}$ with 
$\omega_t^{\nu_0,R,0} := \omega_0^0$
and for $n \geq 1$, $\omega_0^{\nu_n,R,n} := \omega_0^n$, 
\begin{equation}\label{itsyst} 
d{\omega}_t^{\nu_n,R,n}= \big(\nu_n\Delta \omega_t^{\nu_n,R,n} - K_R(\omega_t^{\nu_{n-1},R,n-1})\mathsterling_{u_t^{\nu_{n-1},R,n-1}} \omega_t^{\nu_n,R,n}\big)dt - \displaystyle %
\sum_{i=1}^{\infty}\mathsterling_i\omega_t^{\nu_n,R,n}\circ dW_t^{i,n}
\end{equation}
where $\nu_n = \frac{1}{n}$ is the viscous parameter ($n>0$) and  $u_t^{\nu_{n-1},R,n-1}=curl^{-1}(\omega_t^{\nu_{n-1},R,n-1})$\footnote
{The operator $curl^{-1}$ is the convolution with the Biot-Savart kernel,   see Remark \ref{biotsavartkernel} for details.}. Also $K_R(\omega_t^{\nu_n,R,n}):=f_R(\|\omega_t^{\nu_n,R,n}\|_{k-1,2})$. The corresponding It\^{o} form of equation \eqref{itsyst} is \footnote{%
The stochastic It\^{o} integral is understood here in the usual sense, see \cite{DaPratoZabczyk} .} 
\begin{equation}\label{itsystito}
d\omega_t^{\nu_n,R,n} = (\nu_n\Delta\omega_t^{\nu_n,R,n}+P_{t}^{n-1,n}(\omega_t^{\nu_n,R,n}))dt-\sum_{i=1}^{\infty} \mathsterling_{i}\omega_t^{\nu_n, R,n}dW_t^{i,n},   
\end{equation}
where $P_t^{n-1,n}(\omega_t^{\nu_n,R,n})$ is defined as  
\begin{equation}\label{ptqt}
P_t^{n-1,n}(\omega_t^{\nu_n,R,n}):= -K_R(\omega_t^{\nu_{n-1},R,n-1})\mathsterling_{u_t^{\nu_{n-1},R,n-1}}\omega_t^{\nu_n,R,n}+\frac{1}{2}\displaystyle\sum_{i=1}^{\infty}\mathsterling_{i}^2\omega_t^{\nu_n,R,n}, \ \ t\ge 0.
\end{equation}
\begin{theorem}\label{psns}If $\omega_0^{\nu_n,R,n} \in C^{\infty}(\mathbb{T}^2)$ is a function with null spatial mean, then the two-dimensional stochastic vorticity equation \eqref{itsystito}
admits a unique  global $(\mathcal F_t)_t$-adapted solution
$\omega^{\nu_n,R,n}=\{\omega_t^{\nu_n,R,n},t\in [0,\infty)\}$ in the space 
$C\big([  0,\infty);C^\infty(\mathbb{T}^2) \big)$.  
\end{theorem}
The proof of this theorem is provided in Section \ref{approx.sequence}.

\begin{proposition}\label{tightnessI} 
The laws of the family of solutions $(\omega^{\nu_n,R,n})_{\nu_n\in [0,1]}$ is relatively compact in the space of probability measures over $D([0,T],L^{2}(\mathbb{T}^2))$ for any $T\geq 0$. 
\end{proposition}
The proof of Proposition \ref{tightnessI} is left for Section \ref{sectiontightness}. \\

\noindent \textbf{\textit{Proof of existence of the solution of equation \eqref{truncatedEulere} }} \\
 Using a diagonal subsequence argument we can deduce
  from  Proposition \ref{tightnessI} and the fact that\\ $\displaystyle\lim_{n\rightarrow \infty} \omega_0^{\nu_n,R,n} 
=\omega_0$ the existence of a subsequence $(\omega^{\nu_{n_j}})_{j}$ with $\displaystyle\lim_{j\rightarrow \infty}\nu_{n_j}=0$, 
which is convergent in distribution over  $D([0,\infty),L^2(\mathbb{T}^2))$. We show that the limit of the corresponding distributions is the distribution of a stochastic process that solves \eqref{truncatedEulere}. This justifies the existence of a weak (probabilistic) solution. By using the Skorokhod representation theorem (see \cite{Billingsley} Section 6, pp. 70), there exists a space $(\tilde{\Omega}, \tilde{\mathcal{F}}, \tilde{\mathbb{P}})$ and a sequence of processes $(\tilde{\omega}^{\nu_n,R,n}, \tilde{u}^{\nu_n,R,n}, (\widetilde{W}^{i,n})_i, n=1, \infty)$ which has the same distribution as that of the original converging subsequence and which converges (when $n \rightarrow \infty$) almost surely to a triplet $(\tilde{\omega}^R, \tilde{u}^R, (\widetilde{W}^i)_i )$ in $D([0,T],L^2(\mathbb{T}^2)\times \mathcal{W}^{1,2}(\mathbb{T}^2)\times \mathbb{R}^{\mathbb N})$. Note that $\omega^{\nu_n,R,n}$ and $\tilde{\omega}^{\nu_n,R,n}$ have the same distribution, so that for any test function $\varphi \in C^{\infty}(\mathbb{T}^2)$ 
we have 
\begin{equation}\label{reff}
\begin{aligned}
\langle \tilde{\omega}_t^{\nu_n,R,n}, \varphi\rangle &= \langle \tilde{\omega}_0^{\nu_n,R,n}, \varphi\rangle + \nu_n\displaystyle\int_0^t \langle\tilde{\omega}_s^{\nu_n,R,n},  \Delta\varphi\rangle ds - \displaystyle\int_0^t K_R(\tilde{\omega}_s^{\nu_{n-1,R,n-1}})\langle \tilde{\omega}_s^{\nu_n,R,n},\mathsterling_{u_t^{\nu_{n-1},R,n-1}} \varphi \rangle ds \\
& + \frac{1}{2}\displaystyle\sum_{i=1}^{\infty}\displaystyle\int_0^t \langle\tilde{\omega}_s^{\nu_n,R,n},  \mathsterling_i^2\varphi\rangle ds - \displaystyle\sum_{i=1}^{\infty}\displaystyle\int_0^t \langle \tilde{\omega}_s^{\nu_n,R,n},  \mathsterling_{i}\varphi\rangle d\tilde{W}_s^{i,n}.
\end{aligned}
\end{equation}
Note that there exists a constant $C=C(T)$
such that  
\begin{equation}\label{C(R,T)}
\displaystyle\sup_{n\geq 1}\tilde{\mathbb{E}}\bigg[\displaystyle\sup_{s\in[0,T]}\|\tilde{\omega}^{\nu_n,R,n}_s\|_{k,2}^4\bigg] \leq C,
\end{equation}
where $\tilde{\mathbb{E}}$ is the expectation with respect to $\tilde{\mathbb{P}}$.  We prove this in Lemma \ref{apriori} for the original sequence, but since $\tilde{\omega}^{\nu_n,R,n}$ satisfies the same SPDE, the same a priori estimates hold for $\tilde{\omega}^{\nu_n,R,n}$. Since the space of continuous functions is a subspace of the space of c\`adl\`ag functions and the Skorokhod topology relativised to the space of continuous functions coincides with the uniform topology, it follows that the sequence  $(\tilde{\omega}^{\nu_n,R,n}, \tilde{u}^{\nu_n,R,n}, (\widetilde{W}^{i,n})_i, n=1, \infty)$ converges (when $n \rightarrow \infty$) $\tilde{\mathbb{P}}$-almost surely to $(\tilde{\omega}^R, \tilde{u}^R, (\widetilde{W}^i)_i )$ also in the uniform norm. It also holds that
\[ 
\lim_{n\rightarrow\infty}\tilde{\mathbb{E}}\left[\int_0^t || \tilde{\omega}_s^{\nu_n,R,n}-\tilde{\omega}_s^R||^2ds\right]=0
\]
and since
\begin{equation}\label{tt}
\begin{aligned} 
\displaystyle\sum_{i=1}^{\infty}\displaystyle\tilde{\mathbb{E}}\left[\int_0^t (\langle \tilde{\omega}_s^{\nu_n,R,n}-\tilde{\omega}_s^R, \mathsterling_i\varphi\rangle )^2ds\right] &\le  \displaystyle\sum_{i=1}^{\infty} \|\mathsterling_{i}\varphi\|_{2}^2 \tilde{\mathbb{E}}\left[\int_0^t || \tilde{\omega}_s^{\nu_n,R,n}-\tilde{\omega}_s^R||^2ds\right] \\
&\leq C\|\varphi\|_{1,2}^2 \tilde{\mathbb{E}}\left[\int_0^t || \tilde{\omega}_s^{\nu_n,R,n}-\tilde{\omega}_s^R||^2ds\right]
\end{aligned} 
\end{equation}
also the limit of the right hand side of \eqref{tt} converges to $0$ (we use here the control 
$\sum_{i=1}^{\infty} \|\mathsterling_{i}\varphi\|_{2}^2 \leq C\|\varphi\|_{1,2}^2$ assumed in  \eqref{xicondb}). 
Now Theorem \ref{Kurtzstochint}
 allows us to conclude that 
the sequence of processes 
\begin{equation*}
\left(\tilde{\omega}^{n},\tilde{W}^{i,n},\int_{0}^{\cdot }\langle \tilde{\omega}_{s}^{n},%
\mathsterling_{i}\varphi\rangle d\tilde{W}^{i,n},i\geq 1\right) _{n}
\end{equation*}%
converges in distribution to 
\begin{equation*}
\left(\tilde{\omega}, \tilde{W}^{i}\text{,}\int_{0}^{\cdot }\langle \tilde{\omega}_{s},%
\mathsterling_{i}\varphi\rangle d\tilde{W}^{i},i\geq 1\right) .
\end{equation*}%
By a similar application of the Skorokhod representation theorem, we can also assume that on $(\tilde{\Omega}, \tilde{\mathcal{F}}, \tilde{\mathbb{P}})$, the above convergence
holds also
$\tilde{\mathbb{P}}$-almost surely (as well as in $L^2(\tilde{\mathbb{P}})$). Let us prove the convergence of the remaining terms in \eqref{reff} :

\begin{itemize}
\item $\tilde{\omega}^{\nu_n,R,n}$ converges $\tilde{\mathbb{P}}$-almost surely to  
$\tilde{\omega}^R$ in $D([0,\infty), L^2(\mathbb{T}^2))$. Since $\varphi$ is bounded, it follows that   
$\langle \tilde{\omega}_t^{\nu_n,R,n}, \varphi\rangle \xrightarrow[n \rightarrow \infty]{} \langle \tilde\omega_t^R, \varphi\rangle$ and $\langle \tilde{\omega}_0^{\nu_n,R,n}, \varphi\rangle \xrightarrow[n \rightarrow \infty]{} \langle \tilde\omega_0^R, \varphi\rangle,$ $\tilde{\mathbb{P}}$-almost surely (as well as in $L^2(\tilde{\mathbb{P}})$), for any $\varphi \in C^{\infty}(\mathbb{T}^2)$.
\item The second term on the right hand side of \eqref{reff} converges to $0$ when $n \rightarrow \infty$ because the integral is uniformly bounded in $L^2(\tilde{\mathbb{P}})$ (again, because of \eqref{C(R,T)}) and $\nu_n \rightarrow 0$ when $n \rightarrow \infty$. 

\item Using the fact that $\tilde{u}^{\nu_{n-1},R,n-1}$ is the convolution between $\tilde\omega^{\nu_{n-1},R,n-1}$ and the Biot-Savart kernel we obtain that  $\tilde{u}^{\nu_n,R,n}$ converges to $\tilde{u}^R$, $\tilde{\mathbb{P}}$-almost surely (as well as in $L^2(\tilde{\mathbb{P}})$). Moreover, one can write
   \begin{equation*} 
        \begin{aligned}
        \displaystyle\int_{0}^{t}\langle \tilde{u}_s^{\nu_{n-1},R,n-1} \cdot \nabla{\tilde{\omega}}_s^{\nu_n,R,n}-\tilde{u}_s^R \cdot \nabla\tilde{\omega}_s^R,\varphi \rangle ds & = \displaystyle\int_{0}^{t}\langle (\tilde{u}_s^{\nu_{n-1},R,n-1} - \tilde{u}_s^R) \cdot \nabla{\tilde{\omega}}_s^{\nu_n,R,n}, \varphi\rangle ds \\
        &-\displaystyle\int_{0}^{t}\langle \tilde{u}_s^R \cdot (\nabla\tilde{\omega}_s^{\nu_n,R,n} - \nabla\tilde{\omega}_s^R),\varphi \rangle ds. \\
        \end{aligned} 
   \end{equation*}
   We have 
        \begin{equation*} 
        \begin{aligned} 
        |\langle (\tilde{u}_s^{\nu_{n-1},R,n-1} - \tilde{u}_s^R) \cdot \nabla{\tilde{\omega}}_s^{\nu_n,R,n}, \varphi\rangle| & = |\langle (\tilde{u}_s^{\nu_{n-1},R,n-1} - \tilde u_s^R) \cdot \tilde{\omega}_s^{\nu_n,R,n}, \nabla\varphi\rangle| \\
        &\leq \|\nabla \varphi \cdot \tilde{\omega}_s^{\nu_{n},R,n}\|_2\|\tilde{u}_s^{\nu_{n-1},R,n-1}-\tilde u_s^R\|_2 \xrightarrow[n\rightarrow\infty]{}0
        \end{aligned} 
        \end{equation*}
   and \begin{equation*}
        \begin{aligned}
        |\langle \tilde{u}_s \cdot (\nabla\tilde{\omega}_s^{\nu_n,R,n} - \nabla \tilde{\omega}_s^R),\varphi\rangle| & = |\langle \tilde{u}_s^R \cdot (\tilde{\omega}_s^{\nu_n,R,n} -\tilde{\omega}_s^R),\nabla\varphi\rangle| \\
        & \leq \|\nabla\varphi \cdot \tilde{u}_s^R\|_2\|\tilde{\omega}_s^{\nu_n,R,n} - \tilde{\omega}_s^R\|_2 \xrightarrow[n\rightarrow\infty]{}0
        \end{aligned}
       \end{equation*}
        for $0 \leq t \leq T$.

 \item Lastly, the integrals coming from the It\^{o} correction term are treated in a similar fashion:
\ \begin{equation*}
 \begin{aligned}
 |\langle \xi_i \cdot \nabla(\xi_i \cdot \nabla\tilde{\omega}_s^{\nu_n,R,n})- \xi_i \cdot \nabla(\xi_i \cdot \nabla\tilde\omega_s^{R}), \varphi\rangle| &= |\langle \xi_i \cdot \nabla\tilde{\omega}_s^{\nu_n,R,n} -\xi_i \cdot \nabla\tilde\omega_s^{R}, \xi_i \cdot \nabla\varphi\rangle| \\
 & = |\langle \tilde \omega_s^{\nu_n,R,n} - \tilde\omega_s^{R}, \xi_i \cdot \nabla(\xi_i \cdot \nabla\varphi)\rangle| \\
 & \leq \|\xi_i \cdot \nabla(\xi_i \cdot \nabla\varphi)\|_2 \|\tilde{\omega}_s^{\nu_n,R,n} - \tilde{\omega}_s^R\|_2 \xrightarrow[n\rightarrow\infty]{}0
 \end{aligned}
 \end{equation*}
 since $\|\xi_i \cdot \nabla(\xi_i \cdot \nabla\varphi)\|_2$ is finite by condition \eqref{xicondb} imposed initially on $(\xi_i)_i$.
\end{itemize}
We have shown so far that  there exists a weak/distributional solution in the sense of  Definition \ref{solutions}. part b. on the space $(\tilde{\Omega}, \tilde{\mathcal{F}}, \tilde{\mathbb{P}})$. However, since $\tilde{\omega}^R$ belongs to the space $\mathcal{W}^{k,2}(\mathbb{T}^2) \hookrightarrow C^{k-m}(\mathbb{T}^2)$ the solution is also strong, again, as a solution on $(\tilde{\Omega}, \tilde{\mathcal{F}}, \tilde{\mathbb{P}})$ (and not on the original space). It follows that  $(\tilde{\omega}, \tilde{u}, (\widetilde{W}^i)_i)$ is a martingale solution of the truncated Euler equation \eqref{truncatedEulere}  in the sense of Definition \ref{solutions} part c. Together with the pathwise uniqueness proved in Section \label{uniquenessEulert} \ref{uniquenessEulert} and using the Yamada-Watanabe theorem for the infinite-dimensional setting (see, for instance, \cite{Rockner}) we conclude the existence of a strong solution of the truncated Euler equation. 
Continuity follows immediately: we first apply the Kolmogorov-Čentsov criterion for the approximating process, to control $\mathbb{E}\left[\|\omega_t^{\nu_n,R,n}-\omega_s^{\nu_n,R,n}\|_{k,2}^4\right]$, and then we pass to the limit using an argument similar to the one from Section \ref{contEuler}, to control $\mathbb{E}\left[\|\omega_t^{R}-\omega_s^{R}\|_{k,2}^4\right]$. Therefore $\omega^R \in C([0,T], \mathcal{W}^{k,2}(\mathbb{T}^2))$\footnote{Note that this implies continuity also for the original global solution i.e. $\omega \in C([0,T], \mathcal{W}^{k,2}(\mathbb{T}^2)).$} and the solution exists in the sense of  Definition \ref{solutions}, part a.
Now using the embedding $\mathcal{W}^{k,2}(\mathbb{T}^2) \hookrightarrow C^{k-m}(\mathbb{T}^2)$ with $\ 2\leq m \leq k$ and $k\geq 4$ we conclude that the solution is classical when $k \geq 4$.

\subsection{Proof of Theorem \ref{mainresult2} \label{contEuler}}

We are finally ready to show continuity with respect to initial conditions. As stated in the theorem, let $\omega$, $\tilde\omega$ be two $C\left( [  0,\infty);\mathcal{W}^{k,2}(\mathbb{T}^2)  \right)  $-solutions of equation \eqref{maineqn} and define $A$ as the process $A_t:=\displaystyle\int_0^t\|\omega_s\|_{k,2}ds$, for any $ t\ge 0$. Let $\omega^R$, $\tilde\omega^R$ be their corresponding truncated versions and also let $(\omega_t^{\nu_n,R,n})_{n\geq 0}$ and  $(\tilde\omega_t^{\nu_n,R,n})_{n\geq 0}$ be, respectively, the corresponding sequences constructed as in Section \ref{existenceEulert} on the same space after the application of the Skorokhod representation theorem. By Fatou's lemma, applied twice, we deduce that   
\begin{equation*} 
\begin{aligned}
\mathbb E [e^{-CA_t}||\omega_t-\tilde\omega_t||_{k,2}^2] & \le 
\mathbb E \big[\displaystyle\liminf_n e^{-CA^n_t}||\omega_t^{\nu_n,R,n}-\tilde\omega_t^{\nu_n,R,n}||_{k,2}^2\big] \\
& \leq  \displaystyle\liminf_n \mathbb{E} [e^{-CA^n_t}||\omega_t^{\nu_n,R,n}-\tilde\omega_t^{\nu_n,R,n}||_{k,2}^2]
\end{aligned}
\end{equation*}
where $A^n$ is the process defined by $A_t^n:=\displaystyle\int_0^t\|\omega_s^{\nu_n,R,n}\|_{k,2}ds$, for any $ t\ge 0$.
Following a similar proof with that of the uniqueness of the Euler equation, one then deduces that  
there exists a positive constant $C$ independent of the two solutions and  independent of $R$ and $n$ such that 
\[
\mathbb{E} [e^{-CA^n_t}||\omega_t^{\nu_n,R,n}-\tilde\omega_t^{\nu_n,R,n}||_{k,2}^2]\le
||\omega_0-\tilde\omega_0||_{k,2}^2.
\]
which gives the result. We emphasize that we use here the fact that the processes $(\omega_t^{\nu_n,R,n})_{n\geq 0}$ and $(\tilde\omega_t^{\nu_n,R,n})_{n\geq 0}$ take values in $\mathcal{W}^{k+2,2}(\mathbb{T}^2)$ as an essential ingredient, a property that was not true for either the solution of the Euler equation or its truncated version.  

{\section{Existence, uniqueness, and continuity of the approximating sequence of solutions}\label{approx.sequence}}

\subsection{Existence and uniqueness of the approximating sequence} 

We show that the sequence $(\omega_t^{\nu_n,R,n})_n$ given by formula (\ref{itsystito}), that is, 
\begin{equation*}
d\omega_t^{\nu_n,R,n} = (\nu_n\Delta\omega_t^{\nu_n,R,n}+P_{t}^{n-1,n}(\omega_t^{\nu_n,R,n}))dt-\sum_{i=1}^{\infty} \mathsterling_{\xi_i}\omega_t^{\nu_n, R,n}dW_t^{i,n}   
\end{equation*}
with $\omega_0^{\nu_n,R,n}=\omega_0^n$, is smooth. The equation above is a particular case of equation $(1.1)-(1.2)$ in Chapter 4, Section 4.1, pp.129 in \cite{Rozovskii}. All assumptions required by Theorem 1 and Theorem 2 in \cite{Rozovskii}, Chapter 4, are fulfilled. 
Therefore there exists a unique solution $\omega_t^{\nu_n,R,n}$ which belongs to the class $L^2([0,T], \mathcal{W}^{k,2}(\mathbb{T}^2)) \cap C([0,T], \mathcal{W}^{k-1,2}(\mathbb{T}^2))$ 
and satisfies equation (\ref{itsystito}) for all $t\in[0,T]$ and for all $\omega$ in $\Omega' \subset \Omega$ with $\mathbb{P}(\Omega') = 1.$ 

Furthermore, since the conditions are fulfilled for all $k\in\mathbb{N}$, using Corollary 3
from pp. 141 in \cite{Rozovskii}, we obtain that $\omega_t^{\nu_n,R,n}$ is $\mathbb{P}$ - a.s. in $C\big([0,T], C^{\infty}(\mathbb{T}^2)\big)$. 
Note that $u_t^{\nu_{n-1},R,n-1} \in C^{\infty}(\mathbb{T}^2)$ for any $n \geq 1$, using the Biot-Savart law and an inductive argument. One has $u_t^{\nu_{n-1},R,n-1} = K \star \omega_t^{\nu_{n-1},R,n-1}$ with $K$ being the Biot-Savart kernel defined in Appendix. The convolution between $K$ and $\omega_t^{\nu_{n-1},R,n-1}$ is commutative, so we have
\begin{equation*}
u_t^{\nu_{n-1},R,n-1}(x) = \displaystyle\int_{\mathbb{T}^2}K(y)\omega_t^{\nu_{n-1},R,n-1}(x-y)dy.
\end{equation*}
Since $\omega_t^{\nu_{n-1}}$ is in $C^{\infty}(\mathbb{T}^2)$ by Corollary 3 (at step $n-1$), and using the fact that $K\in L^1(\mathbb{T}^2)$, we conclude that $u_t^{\nu_{n-1},R,n-1}\in C^{\infty}(\mathbb{T}^2)$. This, together with the initial assumption \eqref{xiassumpt}, ensures that all the coefficients of equation \eqref{itsystito} are infinitely differentiable. The uniform boundedness is ensured by the truncation $K_R(\omega_t^{\nu_{n-1},R,n-1})$, as proven in Lemma \ref{apriori} from Appendix.

\subsection{Continuity of the approximating sequence}

\label{solutioncontinuity}
\begin{proposition}\label{cont}
There exists a constant $C=C(T)$ independent of $n$ and $R$ such that    
\begin{equation*}
 \mathbb{E}[\|\omega_t^{\nu_n,R,n} -
\omega_s^{\nu_n,R,n}\|_{L^2}^4] \leq C(t-s)^2,\ \ t,s\in[0,T].
\end{equation*}
In particular, by the Kolmogorov-\v{C}entsov
criterion (see \cite{KaratzasShreve}), the processes $\omega^{\nu_n,R,n}$ have continuous trajectories in $L^2(\mathbb{T}^2)$. 
\end{proposition}

\noindent\textbf{\textit{Proof}} Consider $s\leq t$.
 Then 
\begin{equation}\label{refff}
\omega_t^{\nu_n,R,n}-\omega_s^{\nu_n,R,n} =\nu_n\displaystyle\int_s^t 
\Delta\omega_p^{\nu_n,R,n}dp 
-\displaystyle\int_s^t 
P_{p}^{n-1,n}(\omega_p^{\nu_n,R,n}))dp-\sum_{i=1}^{\infty} \displaystyle\int_s^t \mathsterling_{i}\omega_p^{\nu_n, R,n}dW_p^{i,n}   
\end{equation}
 We will estimate the expected
value of each of these terms. For the first term we have  
\begin{equation*}
\begin{aligned} \mathbb{E}\bigg[\bigg\|\displaystyle\int_s^t\nu_n
\Delta\omega_p^{\nu_n,R,n}dp \bigg\|_{L^2}^4\bigg]&\leq
(t-s)^4\mathbb{E}[\displaystyle\sup_{p\in[s,t]}\|\Delta\omega_p^{\nu_n,R,n}\|_2^4]
\\ & \leq
T^{2}(t-s)^2\mathbb{E}[\displaystyle\sup_{p\in[0,t]}\|\omega_p^{\nu_n,R,n}\|_{k,2}^4]\\ &
\leq C(t-s)^2. \end{aligned}
\end{equation*}Next we have
\begin{equation}
\begin{aligned} &\mathbb{E}\bigg[ \bigg\|\displaystyle\int_s^t K_R(\omega_p^{\nu_{n-1},R,n-1})\mathsterling_{u_t^{\nu_{n-1},R,n-1}} \omega_p^{\nu_n,R,n} dp \bigg\|_{L^2}^4\bigg] \\
&\leq
(t-s)^4\mathbb{E}\bigg[\sup_{p\in[s,t]}\|K_R(\omega_p^{\nu_{n-1},R,n-1})u_p^{\nu_{n-1},R,n-1} \cdot
\nabla\omega_p^{\nu_n,R,n}\|_{L^2}^4\bigg] \\ &\leq T^2(t-s)^2\mathbb{E}\bigg[
\sup_{p\in[0,T]}\|\omega_p^{\nu_n,R,n}\|_{k,2}^4 \bigg] \\ &\leq C(t-s)^2. \end{aligned}\label{zdig}
\end{equation}
\noindent The penultimate inequality is true given that 
\begin{equation*}
\|K_R(\omega_p^{\nu_{n-1},R,n-1})u_p^{\nu_{n-1},R,n-1} \cdot \nabla\omega_p^{\nu_n,R,n}\|_{L^2}^2 \leq
C\|u_p^{\nu_{n-1},R,n-1}\|_{\infty}^2\|\nabla\omega_p^{\nu_n,R,n}\|_2^2 \\
 \leq  C\|\omega_p^{\nu_n,R,n}\|_{k,2}^2 
\end{equation*}
since 
\begin{equation*}
\|\nabla u_p^{\nu_{n-1},R,n-1}\|_{\infty} \leq C\|\nabla u_p^{\nu_{n-1},R,n-1}\|_{k,2}\leq C
\|u_p^{\nu_{n-1},R,n-1}\|_{k+1,2}\leq C\|\omega_p^{\nu_{n-1},R,n-1}\|_{k,2},
\end{equation*}

\noindent and the last term is finite in expectation due to the a priori estimates proved in Lemma \ref{apriori} vi. Similarly, we can prove that
 \[
 \mathbb{E}\bigg[ \bigg\|\displaystyle\int_s^t \frac{1}{2}\displaystyle\sum_{i=1}^{\infty}\mathsterling_{i}^2\omega_p^{\nu_n,R,n} dp \bigg\|_{L^2}^4\bigg] \leq C(t-s)^2. \]
which, together with \eqref{zdig} gives a control on the second term of \eqref{refff}.
For the last term we use the Burkholder-Davis-Gundy inequality and obtain
\begin{equation*}
\begin{aligned}
\mathbb{E}\bigg[\bigg\|\displaystyle\int_s^t\displaystyle\sum_{i=1}^{\infty}
\xi_i \cdot \nabla\omega_p^{\nu_n,R,n} dW_p^{i,n}\bigg\|_{L^2}^4\bigg] &\leq
C\mathbb{E}\bigg[\bigg(\displaystyle\int_s^t \displaystyle\sum_{i=1}^{\infty}
\|\xi_i \cdot \nabla\omega_p^{\nu_n,R,n}\|_2^2 dp\bigg)^2\bigg] \\
& \leq
C(t-s)^2\mathbb{E}\big[\displaystyle\sup_p\displaystyle\sum_{i=1}^{\infty}%
\big\|\xi_i \cdot \nabla\omega_p^{\nu_n,R,n}\big\|_2^4\big] \\
 & \leq C(t-s)^2 \mathbb{E}\big[\displaystyle\sup_p\|\omega_p^{\nu_{n,R,n}}\|_{k,2}^4\big] \\
 & \leq C(t-s)^2
\end{aligned}
\end{equation*}
due to the initial assumption (\ref{xiassumpt})
and the a priori estimates (\ref{apriori}). The conclusion now follows by a direct application of the Kolmogorov-\v{C}entsov criterion.

\section{Relative compactness of the approximating sequence of solutions} \label{sectiontightness}
In this section we prove that the approximating sequence of solutions constructed in Section \ref{existenceEulert} is relatively compact in the space $D([0,T], L^2(\mathbb{T}^2))$\footnote{In fact the paths are continuous in
$L^2(\mathbb{T}^2)$, however Kurtz's criterion only requires càdlàg paths.}. \\
\noindent \textbf{\textit{Proof of Proposition \ref{tightnessI}}} \\
 In order to prove relative compactness we use Kurtz' criterion for relative compactness. For completeness we state the result in Appendix, see Theorem \ref{kurtz}. To do so we need to show that, for every $\eta>0$ there exists a compact set $K_{\eta, t} \subset L^2(\mathbb{T}^2)$ such
that $\displaystyle\sup_{n}\mathbb{P}\big(\omega_t^{\nu_n,R,n}\notin
K_{\eta, t}\big) \leq \eta.$ The compact we use is 
\[
K_{\eta, t} :=\left\{ \omega \in\mathcal{W}^{k,2}(\mathbb{T}^2) | \ \  \|\omega\|_{k,2} <\left(\frac{C}{\eta}\right)^{\frac{1}{4}} \right\}
\] 
where $C$ is the constant appearing in the a priori estimates  (\ref{apriori}).
 By a Sobolev compact embedding theorem,  $K_{\eta, t}$ is a compact set in
$L^2(\mathbb{T}^2)$ and 

\begin{equation*}
\displaystyle\sup_{n}\mathbb{P}\big(\omega_t^{\nu_n,R,n}\notin
K_{\eta,t}\big) = \displaystyle\sup_{n}\mathbb{P}\left( \|\omega_t^{\nu_n,R,n}\|_{k,2} \geq  \left(\frac{C}{\eta}\right)^{\frac{1}{4}} \right)  \leq \displaystyle\sup_{n}\frac{\eta}{C}\mathbb{E}\left[\sup_{t\in [0,T]}\|\omega_t^{\nu_n,R,n}\|_{k,2}^4
\right]\leq \eta.  
\end{equation*}

To prove relative compactness, we need to justify part b)\ of Kurtz' criterion, as per Theorem \ref{kurtz}. For this we will  show that there exists a  family $(\gamma_{\delta}^n)_{0<\delta<1}$ of nonnegative random
variables such that 
\[
\mathbb{E} \big[\|\omega_{t+l}^{\nu_n,R,n} - \omega_t^{\nu_n,R,n}\|_{2}^2   | \mathcal{F}_t^{n}\big] \leq \mathbb{E}\big[%
\gamma_{\delta}^{n} | \mathcal{F}_t^{n} \big]
\] 
with $0 \leq \ell \leq \delta$ and $\displaystyle \lim_{\delta \rightarrow 0}
\sup_{n} \mathbb{E}\big[ \gamma_{\delta}^{n} \big] =0 $ for $t \in
[0, T]$. The filtration $(\mathcal{F}_t^n)_{t}$ corresponds here to the natural filtration $(\mathcal{F}_t^{\omega^{\nu_n,R,n}})_t$. We will use the mild form of equation (\ref{itsystito}), that is
\begin{equation*}
\omega_t^{\nu_n,R,n}=S^n(t)\omega_0^{\nu_n,R,n} -
\displaystyle\int_{0}^{t}S^n(t-s)P_{s}^{n-1,n}(\omega_s^{\nu_n,R,n}))ds-\sum_{i=1}^{\infty} \int_{0}^{t}S^n(t-s)\mathsterling_{i}\omega_s^{\nu_n, R,n}dW_s^{i,n},
\end{equation*}
with $P_s^{n-1,n}$ as defined in (\ref{ptqt}) and $S^n(t):=e^{\nu_n\Delta t}$.
One has 
\begin{equation}\label{mildsum}
\begin{aligned} &\|\omega_{t+l}^{\nu_n,R,n} - \omega_t^{\nu_n,R,n}\|_{2}^2  \leq
C\bigg( \|(S^n(t+l)-S^n(t))\omega_0^{\nu_n,R,n}\|_{2}^2 \\ & +
\bigg\|\displaystyle\int_0^t(S^n(t+l-s)-S^n(t-s))P_s^{n-1,n}(\omega_s^{\nu_n,R,n})ds\bigg\|_{2}^2 +
\bigg\|\displaystyle\int_t^{t+l}S^n(t+l-s)P_s^{n-1,n}(\omega_s^{\nu_n,R,n})ds\bigg\|_{2}^2 \\ & +
\bigg\|\displaystyle\sum_{i=1}^{\infty} \int_{0}^{t}(S^n(t+l-s)-S^n(t-s))\mathsterling_{\xi_i}\omega_s^{\nu_n, R,n}dW_s^{i,n}\bigg\|_{2}^2 +
\bigg\|\displaystyle\sum_{i=1}^{\infty} \int_{t}^{t+l}(S^n(t+l-s)\mathsterling_{i}\omega_s^{\nu_n, R,n}dW_s^{i,n}\bigg\|_{2}^2 \bigg)
\end{aligned}
\end{equation}
We
will estimate each term separately. For the first term we have 
\begin{equation*}
\begin{aligned}
\mathbb{E}\big[\|(S^n(t+l)-S^n(t))\omega_0^{\nu_n,R,n}\|_{2}^2|\mathcal{F}_t^n\big] & \le \|(S^n(l)-1)\omega_0^{\nu_n,R,n}\|_{2}^2 \end{aligned}
\end{equation*}
For the second term, 
\begin{equation*} 
\begin{aligned}
&\mathbb{E}\bigg[\bigg\|\displaystyle\int_0^t(S^n(t+l-s)-S^n(t-s))P_s^{n-1,n}(\omega_s^{\nu_n,R,n})ds\bigg%
\|_{2}^2 \bigg| \mathcal{F}_t^n\bigg] \\
& \leq T\mathbb{E}\bigg[ \displaystyle\int_0^T
\|(S^n(l)-1)P_s^{n-1,n}(\omega_s^{\nu_n,R,n})\|_{2}^2 ds \bigg| \mathcal{F}_t^n\bigg] \end{aligned}
\end{equation*}
For the third term we have 
\begin{equation*}
\begin{aligned}
\mathbb{E}\bigg[\bigg\|\displaystyle\int_t^{t+l}S^n(t+l-s)P_s^{n-1,n}(\omega_s^{\nu_n,R,n})ds\bigg\|_{2}^2
\bigg| \mathcal{F}_t^n\bigg]& \leq
\mathbb{E}\bigg[l^2\displaystyle\sup_{s\in[t,t+l]}\|S^n(t+l-s)P_s^{n-1,n}(\omega_s^{\nu_n,R,n})\|_{2}^2
\bigg| \mathcal{F}_t^n\bigg]\\ & \leq \mathbb{E}\bigg[C
l^2\displaystyle\sup_{s\in[t,t+l]} \|P_s^{n-1,n}(\omega_s^{\nu_n,R,n})\|_{2}^2 \bigg|
\mathcal{F}_t^n\bigg]\\ & \leq \mathbb{E}\bigg[Cl^2\displaystyle\sup_{s \in
[0, T+1]}\|P_s^{n-1,n}(\omega_s^{\nu_n,R,n})\|_{2}^2 \bigg| \mathcal{F}_t^n\bigg] \end{aligned}
\end{equation*}
Aiming to construct the family $(\gamma_{\delta}^n)_{0<\delta<1}$ of nonnegative random variables, we still need to control the two stochastic terms. The first one is more delicate and in order to obtain a suitable control on it we will use the so-called \textit{factorisation formula} (see e.g.  \cite{DaPratoZabczyk} Section 5.3.1). More precisely, we use
the fact that
\begin{equation*}
\displaystyle\int_{a_1}^{a_2}(a_2-r)^{\alpha -1}(r-a_1)^{-\alpha}dr = C({\alpha})
\end{equation*}
where $C(\alpha)$ is a constant which depends on $\alpha >0$ only. Using also the semigroup property $S^n(t-s) = S^n(t-r)S^n(r-s)$ for $s<r<t$, one can
write 
\begin{equation*}
\begin{aligned} \displaystyle\int_0^{t}S^n(t-s)&\mathsterling_i\omega_s^{\nu_n,R,n}dW_s^{i,n} =
C(\alpha)^{-1}\displaystyle\int_0^{t}S^n(t-s)\bigg(\displaystyle\int_s^t(t-r)^{%
\alpha -1}(r-s)^{-\alpha}dr\bigg)\mathsterling_i\omega_s^{\nu_n,R,n}dW_s^{i,n} \\ & = C(\alpha)^{-1}
\displaystyle\int_0^t\bigg(\displaystyle\int_{s}^{t}
(t-r)^{\alpha-1}(r-s)^{-\alpha} S^n(t-r)S^n(r-s)dr\bigg)\mathsterling_i\omega_s^{\nu_n,R,n}dW_s^{i,n} \\ & = C(\alpha)^{-1}
\displaystyle\int_0^t\bigg(\displaystyle\int_{0}^{r}
(t-r)^{\alpha-1}(r-s)^{-\alpha} S^n(t-r)S^n(r-s)\mathsterling_i\omega_s^{\nu_n,R,n}dW_s^{i,n}\bigg)dr \\ & = C(\alpha)^{-1}
\displaystyle\int_0^t\bigg(S^n(t-r)\bigg(\displaystyle\int_0^r(r-s)^{-%
\alpha}S^n(r-s)\mathsterling_i\omega_s^{\nu_n,R,n}dW_s^{i,n}\big)(t-r)^{\alpha-1}\bigg)dr \\ & = C(\alpha)^{-1}
\displaystyle\int_0^tS^n(t-r)z(r)(t-r)^{\alpha -1}dr \end{aligned}
\end{equation*}
where 
\begin{equation*}
z(r):= \displaystyle\int_0^r (r-s)^{-\alpha}S^n(r-s)\mathsterling_i\omega_s^{\nu_n,R,n}dW_s^{i,n}.
\end{equation*}
We choose $\alpha \in (0, 1/2)$ such that all integrals are well-defined. Now the fourth term in \eqref{mildsum} can be estimated as follows: 
\begin{equation*}
\begin{aligned}
&\bigg\|\displaystyle\int_0^t(S^n(t+l-s)-S^n(t-s))\mathsterling_i\omega_s^{\nu_n,R,n}dW_s^{i,n}\bigg\|_{2}^2 \\ &\leq
\displaystyle\int_0^{t}\|(t+l-r)^{\alpha -1}S^n(t+l-r)z(r) - (t-r)^{\alpha
-1}S^n(t-r)z(r)\|_{2}^2dr \\ & = \displaystyle\int_0^t
\|\big((t+l-r)^{\alpha-1}S^n(l)-(t-r)^{\alpha-1}\big)S^n(t-r)z(r)\|_{2}^2dr \\
& \leq \displaystyle\int_0^T
\|\big((t+l-r)^{\alpha-1}S^n(l)-(t-r)^{\alpha-1}\big)S^n(t-r)z(r)\|_{2}^2dr
\end{aligned}
\end{equation*}
and therefore 
\begin{equation*}
\begin{aligned} 
&\mathbb{E}\bigg[ \bigg\|\displaystyle\int_0^t(S^n(t+l-s)-S^n(t-s))\mathsterling_i\omega_s^{\nu_n,R,n}dW_s^{i,n}\bigg\|%
_{2}^2\bigg| \mathcal{F}_t^n\bigg] \\
&\leq \mathbb{E}\bigg[ \displaystyle%
\int_0^T \|\big((t+l-r)^{\alpha-1}S^n(l)-(t-r)^{\alpha-1}\big)%
S^n(t-r)z(r)\|_{2}^2dr \bigg| \mathcal{F}_t^n\bigg].
\end{aligned} 
\end{equation*}
For the
fifth term in \eqref{mildsum} we have  
\begin{equation*}
\begin{aligned} \displaystyle\int_t^{t+l} \|S^n(t+l-s)\mathsterling_i\omega_s^{\nu_n,R,n}\|_{2}^2 ds & \leq
\displaystyle\int_t^{t+l}\|\mathsterling_i\omega_s^{\nu_n,R,n}\|_{2}^2ds \leq l^2 \displaystyle\sup_{s
\in[t, t+l]} \|\mathsterling_i\omega_s^{\nu_n,R,n}\|_{2}^2 \\
&\leq l^2 \displaystyle\sup_{s \in[0, T+1]}
\|\mathsterling_i\omega_s^{\nu_n,R,n}\|_{2}^2, \end{aligned}
\end{equation*}
so 
\begin{equation*}
\mathbb{E}\bigg[\bigg\|\displaystyle\int_t^{t+l}S^n(t+l-s)\mathsterling_i\omega_s^{\nu_n,R,n}dW_s^{i,n} \bigg\|%
_{2}^2\bigg| \mathcal{F}_t^n\bigg] \leq \mathbb{E}\bigg[ l^2 \displaystyle%
\sup_{s \in[0, T+1]} \|\mathsterling_i\omega_s^{\nu_n,R,n}\|_{2}^2\bigg| \mathcal{F}_t^n\bigg].
\end{equation*}

\noindent We can now define 
\begin{equation*}
\begin{aligned} \gamma_l^{\nu_n} &:= \|(S^n(l)-1)\omega_0^{\nu_n,R,n}\|_{2} +
\displaystyle\int_0^T \|(S^n(l)-1)P_s^{n-1,n}(\omega_s^{\nu_n,R,n})\|_{2}^2 ds+ Cl^2\displaystyle\sup_{s
\in [0, T+1]}\|P_s^{n-1,n}(\omega_s^{\nu_n,R,n})\|_{2}^2 \\ &+ l^2 \displaystyle\sup_{s \in[0, T+1]}
\|\mathsterling_i\omega_s^{\nu_n,R,n}\|_{2}^2 + \displaystyle\int_0^T
\|\big((t+l-r)^{\alpha-1}S^n(l)-(t-r)^{\alpha-1}\big)S^n(t-r)z(r)\|_{2}^2dr
\end{aligned}
\end{equation*}
and $\gamma_{\delta}^{n} := \displaystyle\sup_{l \in [0,
\delta]}\gamma_l^{\nu_n}. $ From Lemma \ref{apriori} we deduce that there exist two constants $c_1$ and $c_2$ such that $\mathbb{E}[%
\displaystyle\sup_s\|P_s^{n-1,n}(\omega_s^{\nu_n,R,n})\|_{2}^2] \leq c_1$
and $\mathbb{E}[\displaystyle\sup_s\|\mathsterling_i\omega_s^{\nu_n,R,n}\|_{2}^2] \leq c_2.$ The integrands in the integrals above converge
pointwise to $0$ when $l \rightarrow 0$ due to the strong continuity of the
semigroup $S^n$. 
At the same time, they are bounded by integrable functions,
therefore the convergence is uniform in space by the dominated convergence
theorem. Then the requirement 
\begin{equation*}
\displaystyle \lim_{\delta \rightarrow 0} \sup_{n} 
\mathbb{E}\big[ \gamma_{\delta}^{n} \big] =0
\end{equation*}
is met. In conclusion all the conditions required by Kurtz' criterion are
fulfilled and therefore $(\omega_t^{\nu_n, R, n})_{\nu_n}$ is relatively compact. 

\section{Recovering the solution of the Euler equation in the Yudovich setting} \label{yudovichsection}

As an application of Theorem \ref{mainresult1}, we prove existence of
the solution of the stochastic Euler equation under the relaxed assumption
that $\omega _{0}\in L^{\infty }\left( \mathbb{T}^{2}\right) $. This is the
so-called Yudovich setting, see e.g. \cite{Yudovichoriginal} or \cite{MajdaBertozzi} Section 8.2.
By doing so, we duplicate the result
in \cite{Brzezniak} without the need to impose assumption (\ref{adas}).
The result is as expected: we show existence of a weak solution $\omega
_{t}\in L^{\infty }\left( \mathbb{T}^{2}\right) $ of equation (\ref{maineqn}%
) re-cast in its weak form (\ref{weakdefinition}). Note that whilst the
definition of a weak solution only requires $\omega _{t}\in L^{2}\left( 
\mathbb{T}^{2}\right) $, we are showing here that equation (\ref{maineqn}) has a
solution in the (smaller) space $L^{\infty }\left( \mathbb{T}^{2}\right) $.

The strategy is quite similar with that employed for proving Theorem \ref%
{mainresult1}. We briefly explain here the main steps without going into
details. For an arbitrary $\omega _{0}\in L^{\infty }(\mathbb{T}^{2})$, let $%
(\omega _{0}^{n})_{n}\in \mathcal{W}^{k,2}(\mathbb{T}^{2})$~be a uniformly
bounded sequence such that $\omega _{0}^{n}$ converges to $\omega _{0}$
almost surely. Consider next the sequence of strong solutions $(\omega
^{n})_{n}\in \mathcal{W}^{k,2}(\mathbb{T}^{2})$ of equation (\ref%
{maineqn}), with the corresponding velocities $(u^{n})_{n}\in \mathcal{W}%
^{k+1,2}(\mathbb{T}^{2}).$ Then $(\omega ^{n})_{n}$ is tight as a sequence
with paths in $C([0,T];W^{-2,2}({\mathbb{T}}^{2}))$ (see Lemma \ref%
{newtightness} below).
We can therefore extract a subsequence (we re-index it if necessary) $%
(\omega ^{n})_{n}$ which converges in distribution to $\omega $ over the
space $C([0,T];W^{-2,2}({\mathbb{T}}^{2}))$. By
appealing to Theorem \ref{Kurtzstochint} (via another reduction to a subsequence) we can
deduce that the sequence of processes 
\begin{equation*}
\left( \omega ^{n},W^{i,n},\int_{0}^{\cdot }\mathsterling_{i} \omega _{s}^{n} dW^{i,n},i\geq 1\right) _{n}
\end{equation*}%
converges in distribution to 
\begin{equation*}
\left( \omega ,W^{i}\text{,}\int_{0}^{\cdot }\mathsterling_{i} \omega _{s}^{n} dW^{i},i\geq 1\right) .
\end{equation*}%
Note that the stochastic integrals that appear above are interpreted as $W^{-2,2}$ - processes. Then, using a Skorokhod representation argument similar to the one in Section %
\ref{existenceEulert}, we deduce that the limiting process $\omega $ indeed
satisfies the It\^{o} version of equation \eqref{maineqn}. We do this by
showing that every term in the equation satisfied by $\omega ^{n}$ converges
to the corresponding term in the equation satisfied by $\omega $. The only
dificult term is the nonlinear term, which we analyse in Lemma \ref%
{convergenceofthenonlinearterm} below. This justifies the existence of a
martingale solution of equation \eqref{maineqn} in $C([0,T];W^{-2,2}({%
\mathbb{T}}^{2}))\cap L^{\infty }([0,T];L^{\infty }({\mathbb{T}}^{2}))$
which together with the uniqueness of the weak solution of \eqref{maineqn}
(we use a similar argument as that for Theorem \ref{mainresult1}) provides the existence of a probabilistically strong solution.

To complete the argument we need to prove that the sequence $\left( \omega
^{n}\right) _{n}$ is tight as a sequence with paths in $C([0,T];W^{-2,2}({%
\mathbb{T}}^{2}))$. The main ingredient in the argument is to show that,
just as in the deterministic case, the vorticity $\omega ^{n}$ is propagated
by inviscid flows, and therefore its $L^{p}$-norm is conserved for any $p\in
(0,\infty ]$. For this we characterize the trajectories of the Lagrangian
fluid particles as the solutions of the following stochastic flow 
\begin{equation}
dX_{t}^{n}\left( x\right) =u_{t}^{n}\left( X_{t}^{n}\left( x\right) \right)
dt+\sum_{i=1}^{\infty }\xi _{i}(X_{t}^{n}\left( x\right) )\circ
dW_{t}^{n,i},~X_{0}^{n}\left( x\right) =x\text{, }~x\in \mathbb{T}^{2}.
\label{stocflow}
\end{equation}%
Recall that, following from Theorem \ref{mainresult1}, $(\omega ^{n})_{n}$ $%
\ $has paths in $\mathcal{W}^{k,2}(\mathbb{T}^{2})$ with the corresponding
velocities $(u^{n})_{n}$ having paths in $\mathcal{W}^{k+1,2}(\mathbb{T}%
^{2}) $. By Theorem 4.6.5 in \cite{Kunita97} the mapping $x\rightarrow
X_{t}^{n}\left( x\right) $ is a $C^{k-1}$-diffeomorphism, for all $n\in 
\mathbb{N}$.\footnote{%
Since $t\rightarrow \Vert u_{t}^{n}\Vert _{k-1,\infty }$ is not uniformly
bounded on $\left[ 0,T\right] $, a localization argument is required here
and the fact that $\mathbb{E}\left[ \sup_{t\in \left[ 0,T\right] }\Vert
u_{t}^{n}\Vert _{k-1,\infty }\right] <\infty .$} The stochastic process $%
X_{t}^{n}\left( x\right) $ models the evolution of the Lagrangian particle
path corresponding to a fluid parcel starting from an arbitrary value $~x\in 
\mathbb{T}^{2}$. Each Lagrangian path evolves according to a mean drift flow
perturbed by a random flow which aims to model the rapid oscillations around
the mean. By using the It\^{o}-Wentzel formula, see e.g. \cite{KaratzasShreve}, page 156, one shows that 
\begin{equation*}
d(\omega _{t}^{n}\left( X_{t}^{n}\left( x\right) \right) )=\nabla \omega
_{t}^{n}\left( X_{t}^{n}\left( x\right) \right) \circ dX_{t}^{n}+d\omega
_{t}^{n}\left( X_{t}^{n}\left( x\right) \right) =0.
\end{equation*}%
It follows that $\omega _{t}^{n}\left( X_{t}^{n}\left( x\right) \right)
=\omega _{0}^{n}\left( x\right) \iff \omega _{t}^{n}\left( x\right) =\omega
_{0}^{n}\left( X_{t}^{-n}\left( x\right) \right) $. That is, just as in the
deterministic case, the vorticity is conserved along the particle
trajectories. In particular, this implies that, pathwise, 
\begin{equation}\label{transportlinfty}
\left\vert \left\vert \omega _{t}^{n}\right\vert \right\vert _{\infty
}=\left\vert \left\vert \omega _{0}^{n}\right\vert \right\vert _{\infty
},~~~t\geq 0.
\end{equation}%
In other words, the vorticity remains uniformly bounded for all times and
all realizations. %
%
%
Next, let $h$\ be any measurable function such that $x\rightarrow h\left(
\omega _{t}^{n}\left( x\right) \right) $ is integrable over the torus. Then%
\begin{equation*}
\int_{\mathbb{T}^{2}}h\left( \omega _{t}^{n}\left( x\right) \right) dx=\int_{%
\mathbb{T}^{2}}h\left( \omega _{0}^{n}\left( X_{t}^{-n}\left( x\right)
\right) \right) dx=\int_{\mathbb{T}^{2}}h\left( \omega _{0}^{n}\left(
x\right) \right) \det \left( J_{t}\left( x\right) \right) dx=\int_{\mathbb{T}%
^{2}}h\left( \omega _{0}^{n}\left( x\right) \right) dx
\end{equation*}%
since the determinant of the Jacobian is zero, the fluid being
incompressible. In particular, 
\begin{equation}
\left\vert \left\vert \omega _{t}^{n}\right\vert \right\vert _{p}=\left\vert
\left\vert \omega _{0}^{n}\right\vert \right\vert _{p},~~~t\geq 0,~~~p\in
(0,\infty ].  \label{lpboundstoc}
\end{equation}%
In addition, following the same arguments as in the deterministic case (for
further details and proofs see, for example, \cite{MajdaBertozzi}, pp. 20-23), one can deduce that the vortex lines move with the solution of the
stochastic Euler flow, and that Kelvin's Conservation of circulation is satisfied.%
\footnote{%
Kelvin's conservation of circulation has been shown in \cite{DrivasHolm}.}

\begin{lemma}
\label{newtightness} The sequence $(\omega ^{n})_{n}$ as defined above is
tight as a sequence with paths in $C([0,T];W^{-2,2}({\mathbb{T}}^{2}))$.%
\end{lemma}

\noindent \textbf{\textit{Proof}} Let $\mathcal{Z}$ be the following space\footnote{%
We thank James-Michael Leahy for pointing out this argument to us.}
\begin{equation*}
\mathcal{Z}=\left\{ f\in L^{\infty }\left( 0,T;L^{2}\left( \mathbb{T}%
^{2}\right) \right) \cap C\left( 0,T;W^{-2,2}\left( \mathbb{T}^{2}\right)
\right) |\ \ \Vert f\Vert _{\mathcal{Z}}<\infty \right\}
\end{equation*}%
where 
\begin{equation*}
\Vert f\Vert _{\mathcal{Z}}^{4}:=\sup_{t\in \left[ 0,T\right] }\left\Vert
f\left( t\right) \right\Vert _{L^{2}\left( \mathbb{T}^{2}\right)
}^{4}+\int_{0}^{T}\int_{0}^{T}\frac{\left\Vert f\left( t\right) -f\left(
s\right) \right\Vert _{W^{-2,2}}^{4}}{\left\vert t-s\right\vert ^{2}}dtds\ \
.
\end{equation*}%
Let $B_{\mathcal{Z}}\left( 0,R\right) \subset \mathcal{Z}$ be the closed
ball of radius $R$ in the $\mathcal{Z}$-norm. Then $B_{\mathcal{Z}}\left(
0,R\right) $ is a compact set in $C\left( 0,T;W^{-2,2}\left( \mathbb{T}%
^{2}\right) \right) $. This follows from the generalized Arzela-Ascoli
theorem see Lemma 1 in \cite{SIM} (use the fact that $L^{2}\left( \mathbb{T}%
^{2}\right) $ is compactly embedded in $W^{-2,2}\left( \mathbb{T}^{2}\right) 
$) combined with Lemma 5 in \cite{SIM}.
The tightness then follows as 
\begin{equation*}
\lim_{R\rightarrow \infty }\sup_{n}P\left( \omega^{n}\in B_{\mathcal{Z}%
}\left( 0,R\right) \right) =\lim_{R\rightarrow \infty }\sup_{n}\frac{E\left[
\Vert \omega ^{n}\Vert _{\mathcal{Z}}^{4}\right] }{R^{4}}=0.
\end{equation*}%
This is true as $\sup_{n}E[\Vert \omega ^{n}_{n}\Vert _{\mathcal{Z}%
}^{4}]<\infty $. To show the later claim use the fact that $%
\sup_{n}\Vert \omega^{n}\Vert _{2}<\infty $ due to the existence of the
stochastic flow (see above) and that there exists a constant $C=C\left(
T\right) $ independent of n such that
\begin{equation*}
E\left[ \left\Vert \omega _{t}^{n}-\omega _{s}^{n}\right\Vert _{W^{-2,2}}^{4}%
\right] \leq C\left\vert t-s\right\vert ^{2}.
\end{equation*}

\begin{lemma}
\label{convergenceofthenonlinearterm}Let $(\omega ^{n})_{n}$ be the
sequence constructed above via the Skorokhod representation, and $\varphi \in 
\mathcal{W}^{1,2}(\mathbb{T}^{2})$. Then 
\begin{equation}
\lim_{n\rightarrow \infty }\int_{0}^{t}\langle \omega _{s}^{n},\mathsterling%
_{u_{s}^{n}}\varphi \rangle ds=\int_{0}^{t}\langle \omega _{s},\mathsterling%
_{u_{s}}\varphi \rangle ds.  \label{non}
\end{equation}
\end{lemma}

\noindent\textbf{\textit{Proof}} From (\ref{lpboundstoc}), the sequence $(\omega ^{n})_{n}$
is uniformly bounded and converges to $\omega $ almost surely in $C\left(
0,T;W^{-2,2}\left( \mathbb{T}^{2}\right) \right) $. Let 
\begin{equation*}
M:=\sup_{n>0}\sup_{t\in \left[ 0,T\right] }\left\vert \left\vert \omega
_{t}^{n}\right\vert \right\vert _{L^{2}\left( \mathbb{T}^{2}\right) }<\infty
.
\end{equation*}
Also $\sup_{t\in \left[ 0,T\right] }\left\vert \left\vert \omega
_{t}\right\vert \right\vert _{L^{2}\left( \mathbb{T}^{2}\right) }\leq M.~$We
first deduce that $(\omega _{t}^{n})_{n}$ converges to $\omega _{t}$ in $%
L_w^{2}\left( \mathbb{T}^{2}\right)$ almost surely. To see this, choose an
arbitrary $\varphi $ and $\varepsilon >0.$ Let $\varphi ^{\varepsilon }\in
W^{2,2}\left( \mathbb{T}^{2}\right) $ such that $\left\vert \left\vert
\varphi -\varphi ^{\varepsilon }\right\vert \right\vert _{L^{2}\left( 
\mathbb{T}^{2}\right) }<\varepsilon /4M$ and $N$ such that $\left\vert
\left( \omega _{t}^{n},\varphi ^{\varepsilon }\right) -\left( \omega
_{t},\varphi ^{\varepsilon }\right) \right\vert <\varepsilon /2$ for all $%
n\geq N$. Then, for all $n\geq N,$we have 
\begin{eqnarray*}
\left\vert \left( \omega _{t}^{n},\varphi \right) -\left( \omega
_{t},\varphi \right) \right\vert &=&\left\vert \left( \omega
_{t}^{n},\varphi \right) -\left( \omega _{t}^{n},\varphi ^{\varepsilon
}\right) +\left( \omega _{t}^{n},\varphi \right) -\left( \omega _{t},\varphi
\right) +\left( \omega _{t}^{n},\varphi ^{\varepsilon }\right) -\left(
\omega _{t},\varphi ^{\varepsilon }\right) \right\vert \\
&\leq &\left( \left\vert \left\vert \omega _{t}^{n}\right\vert \right\vert
_{L^{2}\left( \mathbb{T}^{2}\right) }+\left\vert \left\vert \omega
_{t}\right\vert \right\vert _{L^{2}\left( \mathbb{T}^{2}\right) }\right)
\left\vert \left\vert \varphi -\varphi ^{\varepsilon }\right\vert
\right\vert _{L^{2}\left( \mathbb{T}^{2}\right) }+\left\vert \left( \omega
_{t}^{n},\varphi ^{\varepsilon }\right) -\left( \omega _{t},\varphi
^{\varepsilon }\right) \right\vert <\varepsilon ,
\end{eqnarray*}%
hence the claim.

Next observe that, since $u^{n}=K\star \omega ^{n}$, where $K$ is the
Biot-Savart kernel defined in (\ref{expbs}). We can deduce that $(u^{n})_{n}$
is uniformly bounded and converges to $u:=K\star \omega $ almost surely.
Moreover it is uniformly continuous (see Corollary 2.18 in Brzezniak \&
Maurelli). By Arzela-Ascoli theorem we deduce that $u$ is continuous and $%
u^{n}$ converges uniformly to $u$ almost surely. More precisely we have,
almost surely, 
\begin{equation*}
\sup_{t\in \left[ 0,T\right] }\left\vert \left\vert
u_{t}^{n}-u_{t}\right\vert \right\vert _{L^{\infty }\left( \mathbb{T}%
^{2}\right) }=0.
\end{equation*}%
Then 
\begin{eqnarray*}
\left\vert \langle \omega _{s}^{n},\mathsterling_{u_{s}^{n}}\varphi \rangle
-\langle \omega _{s},\mathsterling_{u_{s}}\varphi \rangle \right\vert
&=&\left\vert \langle \omega _{s}^{n},\mathsterling_{u_{s}^{n}}\varphi -%
\mathsterling_{u_{s}}\varphi \rangle +\langle \omega _{s}^{n}-\omega _{s},%
\mathsterling_{u_{s}}\varphi \rangle \right\vert \\
&\leq &\left\vert \langle \omega _{s}^{n},\left( u_{s}^{n}-u_{s}\right)
^{1}\partial _{1}\varphi +\left( u_{s}^{n}-u_{s}\right) ^{2}\partial
_{2}\varphi \rangle \right\vert +\left\vert \langle \omega _{s}^{n}-\omega
_{s},\mathsterling_{u_{s}}\varphi \rangle \right\vert \\
&\leq &\Vert u^{n}-u\Vert _{\infty }\left\vert \left\vert \omega
_{s}^{n}\right\vert \right\vert _{2}\left\vert \left\vert \varphi
\right\vert \right\vert _{1,2}+\left\vert \langle \omega _{s}^{n}-\omega
_{s},\mathsterling_{u_{s}}\varphi \rangle \right\vert
\end{eqnarray*}%
The first term converges to 0 as \ $u^{n}$ converges uniformly to $u$. The
second term converges to 0 as $(\omega ^{n})_n$ is uniformly bounded and converges
to $\omega $ in $L_w^{2}\left( \mathbb{T}^{2}\right)$. The limit in (\ref%
{non}) then follows by the bounded convergence theorem.

\section{Appendix}
In this Appendix we prove the a priori estimates used in the proof of existence of a solution for the Euler equation and we also review some fundamental results mentioned before. We start by introducing the Biot-Savart operator which establishes the connection between the velocity vector field $u$ and the vorticity vector field $\omega$. 
\begin{remark}[The Biot-Savart kernel]\label{biotsavartkernel} The vorticity field
corresponding to a 2D incompressible fluid is conventionally regarded as a
scalar quantity $\omega = curl \ u = \partial_2u_1 - \partial_1u_2$
(formally it is a vector $(0,0, \partial_2u_1 - \partial_1u_2)$ orthogonal
to $u = (u_1, u_2, 0)$ \cite{MajdaBertozzi}). It is known (see \cite%
{MajdaBertozzi}, \cite{Flandolivortex}) that if $\psi : \mathbb{T}^2 \times
[0, \infty) \rightarrow \mathbb{R}$ is a solution for $\Delta\psi = -\omega$
then $u = \nabla^{\perp}\psi$ solves $\omega = curl \ u$, so $u =
-\nabla^{\perp}\Delta^{-1} \omega$. It is worth mentioning that the
existence of a (unique, up to an additive constant) stream function $\psi$ -
and therefore the reconstruction of $u$ from $\omega$ is ensured by the
incompressibility condition $div \ u = 0$ \cite{MajdaBertozzi}. A periodic,
distributional solution of $\Delta \psi = -\omega$ is given by (\cite%
{Flandolivortex}) 
\begin{equation*}
\psi(x) = (G \star \omega)(x)
\end{equation*}
where $G$ is the Green function of the operator $-\Delta$ on $\mathbb{T}^2$, 
$G(x) = \displaystyle\sum_{k \in \mathbb{Z}^2\setminus \{0\}} \frac{e^{ik
\cdot x}}{\|k\|^2}.$ Then the vector field $u=\nabla^{\perp}\psi$ is
uniquely derived from $\omega$ as follows: 
\begin{equation}\label{expbs}
u(x) = (K \star \omega)(x) = \displaystyle\int_{\mathbb{T}^2}
K(x-y)\omega(y)dy
\end{equation}
where $K$ is the so-called \textit{Biot-Savart kernel} 
\begin{equation*}
K(x) = \nabla^{\perp} G(x) = \displaystyle\sum_{k \in \mathbb{Z}^2\setminus
\{0\}} \frac{ik^{\perp}}{\|k\|^2}e^{ik \cdot x}
\end{equation*}
with $k=(k_1, k_2)$, $k^{\perp} = (k_2, -k_1)$. It is known that $G$ is smooth everywhere except at $x=0$, and that $K
\in L^1(\mathbb{T}^2)$. \\
\end{remark}

For the following result we recall a few elementary results of Fourier analysis. We embed $L^{2}\left(  \mathbb{T}^{2}\right)$ into  $L^{2}\left(  \mathbb{T}^{2};\mathbb{C}\right)  $ and consider the basis of functions $\left\{
e^{2\pi i\xi\cdot x};\xi\in\mathbb{Z}^{2}\right\}$. Then every $f\in
L^{2}\left(  \mathbb{T}^{2};\mathbb{C}\right)  $ can be expressed as
\[
f(x)=\sum_{\xi\in\mathbb{Z}^{2}}\widehat{f}\left(  \xi\right)e^{2\pi i \xi\cdot x}
\]
 where  $\widehat{f}\left(  \xi\right)  =\int_{\mathbb{T}^{3}}e^{-2\pi
i\xi\cdot x}f\left(  x\right)  dx$, $\xi\in\mathbb{Z}^{2}$ are the corresponding Fourier coefficients. We have the classical Parseval identity (see e.g. \cite{MajdaBertozzi})
\[
\int_{\mathbb{T}^{2}}\left\vert f\left(  x\right)  \right\vert
^{2}dx=\sum_{\xi\in\mathbb{Z}^{2}}\left\vert \widehat{f}\left(  \xi\right)
\right\vert ^{2}.
\]
 If $v\in L^{2}\left(  \mathbb{T}^{2};\mathbb{R}^{2}\right)
$ is a vector field with components $v_{i}$, $i=1,2$, we write
$\widehat{v}\left(  \xi\right)  =\int_{\mathbb{T}^{2}}e^{2\pi i\xi\cdot
x}v\left(  x\right)  dx$ and we have, in a similar way, that
 $\int_{\mathbb{T}^{2}}\left\vert v\left(  x\right)  \right\vert
^{2}dx=\sum_{\xi\in\mathbb{Z}^{2}}\left\vert \widehat{v}\left(  \xi\right)
\right\vert ^{2}$. Since $u$ and $\omega$ are partial derivatives of other
functions on the torus, they must have zero average:  
\[
\int_{\mathbb{T}^{2}}u^{1}\left(  x\right)  dx=\int_{\mathbb{T}^{2}}u^{2}\left(  x\right)  dx=\int_{\mathbb{T}^{2}}\omega\left(  x\right)  dx=0.
 \]
Hence  
$\widehat{u^{1}}\left(  0,0\right)=\widehat{u^{2}}\left(  0,0\right) 
=\widehat{\omega}\left(  0,0\right)=0$ and the term corresponding to   $\xi=(0,0)$ does not appear in
the Fourier expansion for $u_1$, $u_2$ and, respectively $\omega$.

For every $s\geq0$, let $W^{s,2}\left(\mathbb{T}^{2};\mathbb{C}\right)$ be the fractional Sobolev space of all functions $f\in L^{2}\left(  \mathbb{T}%
^{2};\mathbb{C}\right)$
such that
\[
\sum_{\xi\in\mathbb{Z}^{2}%
}\left\vert \xi\right\vert ^{2s}\left\vert
\widehat{f}\left(  \xi\right)  \right\vert ^{2}<\infty.
\]
It is a simple exercise to show that there exist a constant $C>1$ such that if $s\in \mathbb N$ the  
\[
C^{-1}\left\Vert f\right\Vert _{s,2}^{2}\le \sum_{\xi\in\mathbb{Z}^{2}%
}(1+\left\vert \xi\right\vert ^{2s})\left\vert
\widehat{f}\left(  \xi\right)  \right\vert ^{2}<C\left\Vert f\right\Vert _{s,2}^{2} 
\]
It follows that this definition coincides with the definition given in Section \ref{preliminaries} for integer $s\in\mathbb{N}$.   Therefore we  can extend the norm  $\left\Vert f\right\Vert _{s,2}$ defined for $s\in \mathbb N$ to arbitrary $s>0$ to be given by 
\[
\left\Vert f\right\Vert _{s,2}^{2}= \sum_{\xi\in\mathbb{Z}^{2}%
}(1+\left\vert \xi\right\vert ^{2s})\left\vert
\widehat{f}\left(  \xi\right)  \right\vert ^{2}
\] 
We denote by $W_{\sigma}^{s,2}\left(  \mathbb{T}%
^{2},\mathbb{R}^{2}\right)  $ the space of all zero mean divergence free
(divergence in the sense of distribution) vector fields $v\in L^{2}\left(
\mathbb{T}^{2};\mathbb{R}^{2}\right)  $ such that all components $v_{i}$,
$i=1,2$ belong to $W^{s,2}\left(  \mathbb{T}^{2};\mathbb{R}\right)  $. For
a vector field $v\in W_{\sigma}^{s,2}\left(  \mathbb{T}^{2},\mathbb{R}%
^{2}\right)  $ the norm $\left\Vert v\right\Vert _{s,2}$ is defined by the
identity $\left\Vert v\right\Vert _{s,2}^{2}=\sum_{i=1}%
^{2}\left\Vert v_{i}\right\Vert _{s,2}^{2}$, where $\left\Vert
v_{i}\right\Vert _{s,2}^{2}$ is defined above. We thus have again
$\left\Vert v\right\Vert _{s,2}^{2}:=\sum_{\xi\in\mathbb{Z}%
^{2}\backslash\left\{  0\right\}  }\left\vert \xi\right\vert ^{2s}(\left\vert
\widehat{v}_1\left(  \xi\right)\right\vert+ \left\vert\widehat{v}_1\left(  \xi\right)  \right\vert )^{2}$. For $f\in W^{s,2}\left(
\mathbb{T}^{2};\mathbb{C}\right)  $, we denote by {$\left(  -\Delta\right)
^{s/2}f$ the function of }$L^{2}\left(  \mathbb{T}^{3};\mathbb{C}\right)  $
with Fourier coefficients $\left\vert \xi\right\vert ^{s}\widehat{f}\left(
\xi\right)  $. {For even integers $s\in\mathbb{N}$ , this definition coincides with the classical definition. Similarly, we write $-\Delta^{-1}f$ \ for the function having
Fourier coefficients $\left\vert \xi\right\vert ^{-2}\widehat{f}\left(
\xi\right)  $. We use the same notations for vector fields, meaning that the
operations are made componentwise.}

The Biot-Savart operator is the reconstruction of a zero mean divergence free vector
field $u$ from a divergence free vector field $\omega$ such that
$\operatorname{curl}u=\omega$. As stated in Remark  \ref{biotsavartkernel} on the 2D torus it is given by
$u=-\operatorname{curl}\Delta^{-1}\omega$. It follows that the  Fourier coefficients of $u$ are given by $\widehat{u}\left(  \xi\right)  =\left\vert \xi\right\vert ^{-2}\xi^{\perp}
\widehat{\omega}\left(  \xi\right)$, where $\xi^{\perp}=(\xi_1,\xi_2)^{\perp}=(\xi_2,-\xi_1)$. 

In the next proposition we highlight the smoothing properties of the Biot-Savart kernel $K$.  

\begin{proposition} [The Biot-Savart law, \cite{CFH}, \cite{MajdaBertozzi}] \label{biotsavartlaw}
Let $u$ be the divergence-free, zero average, vector field defined as  $u = -curl
\Delta^{-1} \omega$. Then, for any $s\ge 0$, there exists a constant $C_{s,2}$, independent of $u$ such that 
\begin{equation} \label{biotsavarteqn}
\|u\|_{s+1, 2} \leq C_{s,2}\|\omega\|_{s,2}.
\end{equation}%
\end{proposition}

\noindent \textbf{\textit{Proof}} Using the definition given above of $\left\Vert u\right\Vert
_{s+1,2}^{2}$, the formula which relates $\widehat{u}\left(
\xi\right)  $ and $\widehat{\omega}\left(  \xi\right)  $ and the rule
$\left\vert a\times b\right\vert \leq\left\vert a\right\vert \left\vert
b\right\vert $, we get
\begin{eqnarray*}
\left\Vert u\right\Vert
_{s+1,2}^{2}&=&\sum_{\xi\in\mathbb{Z}%
^{2}\backslash\left\{  0\right\}  }(1+\left\vert \xi\right\vert ^{2s+2})\left\vert
\widehat{u}\left(  \xi\right)  \right\vert ^{2}\\
&\leq &\sum_{\xi\in\mathbb{Z}%
^{2}\backslash\left\{  0\right\}  }(1+\left\vert \xi\right\vert ^{2s+2})\left\vert
\xi\right\vert ^{-4}\vert \xi^\perp\widehat{\omega}\left(  \xi\right)
\vert ^{2}\\
&\leq&\sum_{\xi\in\mathbb{Z}^{2}\backslash\left\{  0\right\}
}(\left\vert \xi\right\vert ^{-2}+\left\vert \xi\right\vert ^{2s})\left\vert
\widehat{\omega}\left(  \xi\right)  \right\vert ^{2}\\
&\leq&\sum_{\xi\in\mathbb{Z}^{2}\backslash\left\{  0\right\}
}(1+\left\vert \xi\right\vert ^{2s})\left\vert
\widehat{\omega}\left(  \xi\right)  \right\vert ^{2}
\end{eqnarray*}
and the latter is precisely equal to $\left\Vert \omega\right\Vert
_{s,2}$.

\begin{remark}\label{normequiv}
The norm $\|\cdot\|_{m,2}$ is equivalent to the norm defined as $||| f |||:= \|f\|_2 + \|D^{m}f\|_2$, therefore it is enough to show that all properties hold for the $L^2$ norm of $f$ and for the $L^2$ norm of the maximal derivative $D^mf$ (see \cite{Brezis} pp. 217).  
\end{remark}

 Let $\omega$ be the solution of the Euler equation \eqref{maineqn} and $\omega^{\nu_n,R,n}$ the solution of the linear approximating equation \eqref{itsyst}.
In the following lemmas we collect a number of identities and a priori 
estimates. Lemma \ref{correctedlemma} in particular in proving that the solution of the Euler equation is global. 

\begin{lemma}\label{correctedlemma}
There exists a constant $\mathcal{C}=\mathcal{C}(\omega_0,T)$ which is independent of the truncation radius $R$ such that
\begin{equation*}
\mathbb{E}\left[ \displaystyle\sup_{s \in [0,t]} \left( \ln (e + \|\omega_s\|_{k,2}^2)\right)\right] \leq \mathcal{C}(\|\omega_0\|_{k,2},T).
\end{equation*}
\end{lemma}

\noindent\textbf{\textit{Proof}}
By the Biot-Savart law
\begin{equation*}
\begin{aligned}
u(x) =  \displaystyle\int_{\mathbb{T}^2} K(x-y)\omega(y)dy = - \frac{1}{4\pi} \displaystyle\int_{\mathbb{T}^2} \frac{(x-y)^{\perp}}{|x-y|^2} \omega(y)dy
\end{aligned}
\end{equation*}
We use the truncation
\begin{equation*}
z_{\epsilon}(x) = \begin{cases}
    1, & \text{if $|x|\leq \epsilon$}\\
    0, & \text{if $|x| >  2\epsilon$}
  \end{cases}
\end{equation*}
 with $|\nabla z_{\epsilon}(x)| \leq \frac{C}{\epsilon}, \epsilon \leq 1$. Let
\begin{equation*}
u^1(x) := -\frac{1}{4 \pi} \displaystyle\int_{\mathbb{T}^2} z_{\rho}(x-y)K(x-y)\omega(y)dy 
\end{equation*} 
\begin{equation*} 
u^2(x) := -\frac{1}{4 \pi} \displaystyle\int_{\mathbb{T}^2}\left(1- z_{\epsilon}(x-y)\right)K(x-y)\omega(y)dy.
\end{equation*}
Then
\begin{equation*}
\begin{aligned} 
\nabla u^1(x) &:= - \frac{1}{4\pi} \displaystyle\int_{\mathbb{T}^2} z_{\epsilon}(x-y)K(x-y)\nabla \omega(y) dy \\
& \leq \left( \displaystyle\int_{\mathbb{T}^2} \left( z_{\epsilon}(x-y)K(x-y)\right)^{4/3}\right)^{3/4} \left(\displaystyle\int_{\mathbb{T}^2} \left( \nabla \omega(y)\right)^4dy\right)^{1/4} \\
& \leq \|K\|_{L^{4/3}(\{y:|x-y| < 2\epsilon\})} \|\nabla \omega\|_{L^{4}(\{y:|x-y| < 2\epsilon\})} \\
& \leq C \epsilon^{2/3} \|\nabla \omega\|_{4} \\
& \leq C \epsilon^{2/3} \|u\|_{3,2}
\end{aligned} 
\end{equation*}
given that $\|\nabla \omega\|_4 \leq C \|\nabla \omega\|_{1,2} \leq C\|\omega\|_{2,2} \leq C \|u\|_{3,2}$ by the Sobolev embedding $W^{1,2} \hookrightarrow L^4$ and the Biot-Savart law. 
So
\begin{equation*}
\|\nabla u^1(x)\|_{\infty} \leq  C \epsilon^{2/3} \|u\|_{3,2}.
\end{equation*}
For the second integral 
\begin{equation*}
\begin{aligned}
\nabla u^2(x) := - \frac{1}{4\pi} \displaystyle\int_{\mathbb{T}^2} \nabla \left( (1 - z_{\epsilon}(x-y))K(x-y)\right) \omega(y)dy
\end{aligned}
\end{equation*}
we use the fact that 
\begin{equation*}
|\nabla K(x-y)| \leq \frac{C}{|x-y|^2}
\end{equation*}
so 
\begin{equation*}
\|\nabla u^2(x)\|_{\infty} \leq C \left(1 - \ln \epsilon \right)\|\omega\|_{\infty}.
\end{equation*} 
That is
\begin{equation*}
\|\nabla u\|_{\infty} \leq C \left(\epsilon^{2/3}\|u\|_{3,2} + (1- \ln \epsilon) \|\omega\|_{\infty}\right). 
\end{equation*}
If $\|u\|_{3,2} \leq 1$ we choose $\epsilon = 1$, otherwise we choose $\epsilon = \frac{1}{\|u\|_{3,2}^{3/2}}$. We have
\begin{equation*}
\|\nabla u\|_{\infty} \leq C \left( 1 + \left(1 + \frac{3}{2} \ln \|u\|_{3,2}\right)\|\omega\|_{\infty}\right)
\end{equation*}
Using the inequality $1 + \ln^{+} x \leq C \ln(e+x) $ for sufficiently large $C$, we have
\begin{equation*}
\begin{aligned}
\|\nabla u_t\|_{\infty} & \leq C \left( 1 + \left(1 + \ln^{+}\|u_t\|_{3,2}\right)\|\omega_t\|_{\infty}\right) \\
&  \leq C \left( 1 + \ln( e + \|u_t\|_{3,2})\|\omega_t\|_{\infty}\right)
\end{aligned} 
\end{equation*}
and using the transport property and Biot-Savart $\|u_t\|_{3,2} \leq \|\omega_t\|_{2,2} \leq \|\omega_t\|_{2,2}$ we have
\begin{equation} \label{graduestim}
\begin{aligned}
\|\nabla u_t\|_{\infty} & \leq C \left( 1 + \ln( e + \|\omega_t\|_{k,2})\|\omega_0\|_{\infty}\right)
\end{aligned} 
\end{equation}
By the It\^{o} formula
\begin{equation*}
\begin{aligned}
d \ln (e + \|\omega_t\|_{k,2}^2) & \leq \frac{1}{e + \|\omega_t\|_{k,2}^2} d\|\omega_t\|_{k,2}^2   - \frac{2}{\left( e + \|\omega_t\|_{k,2}^2\right)^2} \displaystyle\sum_{k=1}^{\infty} \left( |\langle \partial^k \mathsterling_i\omega_t, \partial^k \omega_t\rangle| + |\langle \mathsterling_{i}\omega_t, \omega_t\rangle| \right)^2 dt \\
& \leq  \frac{C}{e + \|\omega_t\|_{k,2}^2} (\|\nabla u_t\|_{\infty} + \|\omega_t\|_{\infty})\|\omega_t\|_{k,2}^2 dt + dY_t \\
& \leq C \left( \|\omega_t\|_{\infty} + C_1 \left( 1 + \|\omega_t\|_{\infty}\ln (e + \|\omega_t\|_{k,2}) \right)\right) dt + dY_t \\
& \leq C \left( C_1 + \|\omega_t\|_{\infty} + \|\omega_t\|_{\infty}\ln (e + \|\omega_t\|_{k,2}^2)\right)dt + dY_t \\
& \leq \left( C_1\ln (e + \|\omega_t\|_{k,2}^2) + \|\omega_t\|_{\infty}\ln (e + \|\omega_t\|_{k,2}^2) + \|\omega_t\|_{\infty}\ln (e + \|\omega_t\|_{k,2}^2)\right)dt + dY_t \\
& \leq C_2 (1 + \|\omega_t\|_{\infty})\ln (e + \|\omega_t\|_{k,2}^2)dt + dY_t \\
\end{aligned} 
\end{equation*}
where
\begin{equation*}
\begin{aligned}
Y_t := \displaystyle\sum_{i=1}^{\infty} \displaystyle\int_0^t \frac{2\left(\langle \partial^k \mathsterling_{i}\omega_s, \partial^k \omega_s\rangle + \langle \mathsterling_{i}\omega_s, \omega_s\rangle \right)}{ e + \|\omega_s\|_{k,2}^2} dW_s^k.
\end{aligned}
\end{equation*}
Define
\begin{equation*}
Z_t : = \displaystyle\int_0^t (1 + \|\omega_s\|_{\infty}) ds.
\end{equation*}
We have
\begin{equation}\label{beforegronwall1}
e^{-C_2 Z_t} \ln (e + \|\omega_t\|_{k,2}^2) \leq \ln (e + \|\omega_0\|_{k,2}^2) + \displaystyle\int_0^t e^{-C_2Z_s} dY_s.
\end{equation} 
By the Burkholder-Davis-Gundy inequality
\begin{equation*}
\mathbb{E}\left[ \displaystyle\sup_{s \in [0,t]} \left| \displaystyle\int_0^s dY_r\right|\right] \leq \alpha \mathbb{E}\left[ \left\langle \displaystyle\int_0^{\cdot} dY_r \right\rangle_t^{1/2} \right]
\end{equation*} 
We control the quadratic variation of the stochastic integral using the fact that 
\begin{equation*}
\displaystyle\sum_{i=1}^{\infty}\|\xi_i\|_{3,2} < \infty 
\end{equation*}
and 
\begin{equation*}
\left| \langle \mathsterling_{i}\omega_t, \omega_t\rangle + \langle \partial^k\mathsterling_{i}\omega_t, \partial^k\omega_t \rangle\right| \leq C\|\xi_i\|_{3,2}^2\|\omega_t\|_{k,2}^2 < \infty
\end{equation*}
since
\begin{equation*}
\left| \langle \mathsterling_{i}\omega_t, \omega_t\rangle \right| \leq \|\nabla \xi_i\|_{\infty}\|\omega_t\|_2^2
\end{equation*}
and
\begin{equation*}
\left| \langle \partial^k\mathsterling_{i}\omega_t, \partial^k\omega_t\rangle\right| \leq \|\xi_i\|_{3,2}\|\omega_t\|_{k,2}^2.
\end{equation*}
We write
\begin{equation*}
\begin{aligned} 
\left[ \displaystyle\int_0^{\cdot} dY_s\right]_{t} &\leq \displaystyle\sum_{i=1}^{\infty} \displaystyle\int_0^t \frac{4\left(\langle \partial^k \mathsterling_{i}\omega_s, \partial^k \omega_s\rangle + \langle \mathsterling_{i}\omega_s, \omega_s\rangle \right)^2}{ \left(e + \|\omega_s\|_{k,2}^2 \right)^2} dW_s^k \\
& \leq 4C \displaystyle\sum_{i=1}^{\infty} \|\xi_i\|_{3,2}^2 \displaystyle\int_0^t \frac{\|\omega_s\|_{k,2}^4}{\left(e + \|\omega_s\|_{k,2}^2 \right)^2}ds \\
& \leq Ct. 
\end{aligned} 
\end{equation*}
In consequence, the Burkholder-Davis-Gundy inequality gives
\begin{equation*}
\mathbb{E}\left[ \displaystyle\sup_{s\in[0,t]} \left| \displaystyle\int_0^s dY_r\right| \right] \leq C \sqrt{t}.
\end{equation*}
Using this in \ref{beforegronwall1} we obtain
\begin{equation*}
\mathbb{E}\left[ \displaystyle\sup_{s \in [0,t]} \left(e^{-C\displaystyle\int_0^s(1+\|\omega_r^R\|_{\infty})dr} \ln (e + \|\omega_s^R\|_{k,2}^2)\right)\right] \leq \ln (e + \|\omega_0\|_{k,2}^2) + C\sqrt{t}.
\end{equation*}
By the Fatou lemma this holds also  for the original solution $\omega$ 
\begin{equation*}
\begin{aligned}
&\mathbb{E}\left[ \displaystyle\sup_{s \in [0,t]} \left(e^{-C\displaystyle\int_0^s(1+\|\omega_r\|_{\infty})dr} \ln (e + \|\omega_s\|_{k,2}^2)\right)\right] \\
&\leq \mathbb{E}\left[ \displaystyle\sup_{s \in [0,t]} \displaystyle\lim\inf_{R \rightarrow \infty}\left(e^{-C\displaystyle\int_0^s(1+\|\omega_r^R\|_{\infty})dr} \ln (e + \|\omega_s^R\|_{k,2}^2)\right)\right] \leq \ln (e + \|\omega_0\|_{k,2}^2) + C\sqrt{t}.\\
\end{aligned}
\end{equation*}
Using the transport property (see relation \eqref{transportlinfty} in Section \ref{yudovichsection}), $\|\omega_t\|_{\infty} = \|\omega_0\|_{\infty}$, 
we conclude that there exists a constant $\mathcal{C}(\|\omega_0\|_{k,2},T)$ which is independent of the truncation radius $R$ such that 
\begin{equation*}
\mathbb{E}\left[ \displaystyle\sup_{s \in [0,t]} \left( \ln (e + \|\omega_s\|_{k,2}^2)\right)\right] \leq \mathcal{C}(\|\omega_0\|_{k,2},T).
\end{equation*}
\begin{lemma}\label{apriori}$\left.\right. $\\ 
The following properties hold: 
\begin{enumerate}
\item[i.] For any $f\in \mathcal{W}^{2,2}(\mathbb{T}^2)$ we have
\[
\big\langle f, \mathsterling_{i}^2
f \big\rangle  + \langle\mathsterling_{i} f, \mathsterling_{i} f \rangle=0.
\]
\item[ii.] If  $\omega_t, \omega_t^{\nu_n,R,n}\in \mathcal{W}^{k,2}(\mathbb{T}^2)$  then the following $L^2$ - transport formulae\footnote{Actually analogous $L^p$ - transport formulae hold, see Section \ref{yudovichsection}.}. hold $\mathbb{P}$ - almost surely:
\begin{equation*}
\|\omega_t\|_{L^2} =\|\omega_0\|_{L^2} \ \ \ \mathrm{and}\ \ \ \ 
\|\omega_t^{\nu_n,R,n}\|_{L^2}\leq
\|\omega_0\|_{L^2}.
\end{equation*}

\item[iii.] If  $\omega\in \mathcal{W}^{k,2}(\mathbb{T}^2)$, then  $(P_t^{n-1,n})_t$ defined in \ref{ptqt} and $(\mathsterling_i\omega_t^{\nu_n,R,n})_t$ are processes with paths taking values in $L^2(\mathbb{T}^2)$. 

\item[iv.] There exists a constant $C_1$ such that: 
\begin{equation*}
\big|\big\langle \partial^k \omega_t^{\nu_n, R, n}, \partial^k \big( \mathsterling_i^2\omega_t^{\nu_n, R, n}\big)\big\rangle + \big\langle\partial^k\big(\mathsterling_i\omega_t^{\nu_n, R, n}\big),\partial^k\big(\mathsterling_i\omega_t^{\nu_n, R, n}%
\big)\big\rangle \big| \leq C_1\|\omega_t^{\nu_n, R, n}\|_{k,2}^2.
\end{equation*}
\item[v.] There exist some constants $C_2$
and $C_2^{\prime}$ such that: 
\begin{equation*}
\big|\big\langle \partial^k \omega_t^{\nu_n, R, n}, K_R(\omega_t^{\nu_{n-1}, R, n-1}) \partial^k\big(\mathsterling_{u_t^{\nu_{n-1}, R, n-1}}
\omega_t^{\nu_n, R, n}\big)\big \rangle \big| \leq
C_2\|\partial^{k}u_t^{\nu_n, R, n-1}\|_2^a\|u_t^{\nu_n, R, n-1}\|_2^{1-a}\|\omega_t^{\nu_n, R, n}%
\|_{k,2}^2
\end{equation*}
with $0<a\leq 1$, and 
\begin{equation*}
|\langle \partial^k \omega_t^{\nu_n, R, n}, K_R(\omega_t^{\nu_{n-1}, R, n-1})\partial^k\big(\mathsterling_{u_t^{\nu_{n-1}, R, n-1}}
\omega_t^{\nu_n, R, n}\big) \rangle| \leq C_2^{\prime }\|\omega_t^{\nu_n, R, n}\|_{k,2}^2.
\end{equation*}%
\item[vi.] For arbitrary $p\ge 2$, there exists a constant $C_p(T)$ independent of $n$ such that 
\begin{equation*}
\mathbb{E}\big[\displaystyle\sup_{t \in [0,T]}\|\omega_t^{\nu_n,R,n}\|_{k,2}^p\big] \leq C_p(T). 
\end{equation*}

\end{enumerate}
\end{lemma}
\noindent\textbf{\textit{Proof of Lemma \ref{apriori}}}\\
i.  Since the dual of $\mathsterling_i$ is $-\mathsterling_i$ by Remark \ref{dual}, observe that 
                     $$
                     \begin{aligned}
                     \big\langle\mathsterling_{i}^2f, f \big \rangle+ \big\langle \mathsterling_{i}f, \mathsterling_{i}f \big\rangle&=
                     \big\langle\mathsterling_{i}f, \mathsterling_{i}^{\star}f \big \rangle+ \big\langle \mathsterling_{i}f, \mathsterling_{i}f\big\rangle \\
                     &=\big\langle\mathsterling_{i}f, -\mathsterling_{i}f\big \rangle+ \big\langle \mathsterling_{i}f, \mathsterling_{i}f\big\rangle=0. 
                     \end{aligned}
                    $$
This is an intrinsic property of the operator $\mathsterling$, which holds even when $f$ is not a solution of the Euler equation or of the approximating sequence. \\
ii. By It\^{o} formula 
\begin{equation*}
\begin{aligned} d\|\omega_t^{\nu_n,R,n}\|_2^2 & = -2 \displaystyle\sum_i\langle
\omega_t^{\nu_n,R,n}, \mathsterling_i\omega_t^{\nu_n,R,n}\rangle dW_t^{i,n} - 2\langle
\omega_t^{\nu_n,R,n}, K_R(\omega_t^{\nu_{n-1},R,n-1})\mathsterling_{u_t^{\nu_{n-1},R,n-1}}\omega_t^{\nu_n,R,n}\rangle dt\\
& + 2\langle
\omega_t^{\nu_n,R,n}, \nu_n\Delta\omega_t^{\nu_n,R,n}\rangle dt + \displaystyle\sum_i
\bigg(\langle \omega_t^{\nu_n,R,n}, \mathsterling_i^2\omega_t^{\nu_n,R,n}\rangle dt + \langle \mathsterling_i\omega_t^{\nu_n,R,n},
\mathsterling_i\omega_t^{\nu_n,R,n}\rangle dt\bigg). \end{aligned}
\end{equation*}
The last sum has been proved to be $0$ at i. and $\langle \omega_t^{\nu_n,R,n},
\mathsterling_{u_t^{\nu_{n-1},R,n-1}}\omega_t^{\nu_n,R,n}\rangle dt = 0$ by integration by parts
since $u_t^{\nu_{n-1},R,n-1}$ is divergence-free. Similarly $\langle \omega_t^{\nu_n,R,n}, \mathsterling_i\omega_t^{\nu_n,R,n}\rangle = 0$ since the vector fields $\xi_i$ are assumed to be divergence-free. 
Therefore 
\begin{equation*}
\|\omega_t^{\nu_n,R,n}\|_2 = \|\omega_0\|_2^2 - 2\displaystyle\int_0^t \nu_n (\nabla\omega_s^{\nu_n,R,n})^2ds \leq \|\omega_0\|_2^2 \ \ \ \ \ \ \mathbb{P}\hbox{ - almost surely}.
\end{equation*}
The calculations are similar for the Euler equation, but there are no viscous terms, hence 
\begin{equation*}
\|\omega_t\|_2 = \|\omega_0\|_2^2 \ \ \ \ \ \ \mathbb{P}\hbox{ - almost surely}.
\end{equation*}
iii. One has 
\begin{equation*}
\|u_t^{\nu_{n-1},R,n-1} \cdot \nabla \omega_t^{\nu_n,R,n}\|_{2}^2 \leq
C\|u_t^{\nu_{n-1},R,n-1}\|_{\infty}^2 \|\nabla\omega_t^{\nu_n,R,n}\|_2^2 \leq
C\|\omega_t^{\nu_n,R,n}\|_{k,2}^4
\end{equation*}
due to H\"{o}lder's inequality and the following inequalities: \newline
a) $\|u_t^{\nu_{n-1},R,n-1}\|_{\infty} \leq C \|\nabla u_t^{\nu_{n-1},R,n-1}\|_{\infty}$ by
Poincar\'e's inequality\\
(or using directly $
\|u_t^{\nu_{n-1},R,n-1}\|_{\infty} \leq \|u_t^{\nu_{n-1},R,n-1}\|_{k,2} \le \| \omega_t^{\nu_{n-1},R,n-1}\|_{k-1,2}
$)

\noindent b) $\|\nabla u_t^{\nu_{n-1},R,n-1}\|_{\infty} \leq C \|\nabla u_t^{\nu_{n-1},R,n-1}\|_{k,2}$ by the
Sobolev embedding theorem $\mathcal{W}^{k,2}(\mathbb{T}^2) \hookrightarrow
L^{\infty}(\mathbb{T}^2)$ (\cite{Adams}, Theorem 4.12, case A, $m=k, p=n=2,
q= \infty$). \newline
c) $\|\nabla u_t^{\nu_{n-1},R,n-1}\|_{k,2} \leq \|u_t^{\nu_{n-1},R,n-1}\|_{k+1,2}$ by the definition
of the Sobolev norm \newline
d) $\|u_t^{\nu_{n-1,R,n-1}}\|_{k+1,2} \leq C\|\omega_t^{\nu_{n-1},R,n-1}\|_{k,2}$ by the Biot-Savart
law.

All the other terms which involve $\xi_i$ stay in $L^2$
according to the initial assumptions \ref{xiassumpt}, therefore the conclusion holds also for the process $(\mathsterling_u\omega_t^{\nu_n,R,n})_t$. \\

\vspace{2mm} 
\noindent iv. Let us denote $\omega_t^{\nu_n,R,n}$ shortly by $\omega_t^{\nu}$. We have to estimate 
\begin{equation*}
\begin{aligned} \langle \partial^k\mathsterling_i^2\omega_t^{\nu},
\partial^k\omega_t^{\nu}\rangle + \langle \partial^k \mathsterling_i\omega_t^{\nu},
\partial^k \mathsterling_i\omega_t^{\nu} \rangle. \end{aligned}
\end{equation*}
Remark that 
\begin{equation*}
\partial^k\mathsterling_i^2\omega_t^{\nu}=\big((\partial^k\mathsterling_i)%
\mathsterling_i\big)\omega_t^{\nu}=\big((\partial^k\mathsterling_i+\mathsterling%
_i\partial^k-\mathsterling_i\partial^k)\mathsterling_i\big)\omega_t^{\nu}=((%
\mathsterling_i\partial^k+L_i)\mathsterling_i)\omega_t^{\nu}
\end{equation*}
where $L_i:=\partial^k\mathsterling_i-\mathsterling_i\partial^k$. One can
write 
\begin{equation*}
\begin{aligned} \langle \partial^k\mathsterling_i^2\omega_t^{\nu},
\partial^k\omega_t^{\nu}\rangle + \langle \partial^k \mathsterling_i\omega_t^{\nu},
\partial^k \mathsterling_i\omega_t^{\nu} \rangle&=\langle
(\mathsterling_i\partial^k\mathsterling_i)\omega_t^{\nu},
\partial^k\omega_t^{\nu}\rangle + \langle(L_i\mathsterling_i)\omega_t^{\nu},
\partial^k\omega_t^{\nu}\rangle + \langle \partial^k \mathsterling_i\omega_t^{\nu},
\partial^k \mathsterling_i\omega_t^{\nu} \rangle\\ &=\langle
\partial^k\mathsterling_i\omega_t^{\nu},
{\mathsterling_i}^{\star}\partial^k\omega_t^{\nu}\rangle +
\langle(L_i\mathsterling_i)\omega_t^{\nu}, \partial^k\omega_t^{\nu}\rangle + \langle
\partial^k \mathsterling_i\omega_t^{\nu}, \partial^k \mathsterling_i\omega_t^{\nu}
\rangle\\ &=\langle \partial^k\mathsterling_i\omega_t^{\nu},
-\mathsterling_i\partial^k\omega_t^{\nu}\rangle +
\langle(L_i\mathsterling_i)\omega_t^{\nu}, \partial^k\omega_t^{\nu}\rangle + \langle
\partial^k \mathsterling_i\omega_t^{\nu}, \partial^k \mathsterling_i\omega_t^{\nu}
\rangle\\
&=\langle\partial^k\mathsterling_i\omega_t^{\nu},-(\partial^k\mathsterling_i-%
\partial^k\mathsterling_i+\mathsterling_i\partial^k) \omega_t^{\nu}\rangle +
\langle(L_i\mathsterling_i)\omega_t^{\nu}, \partial^k\omega_t^{\nu}\rangle + \\
&+\langle \partial^k \mathsterling_i\omega_t^{\nu}, \partial^k
\mathsterling_i\omega_t^{\nu} \rangle \\
&=-\langle\partial^k\mathsterling_i\omega_t^{\nu},\partial^k\mathsterling_i%
\omega_t^{\nu}\rangle +
\langle\partial^k\mathsterling_i\omega_t^{\nu},L_i\omega_t\rangle +
\langle(L_i\mathsterling_i)\omega_t^{\nu}, \partial^k\omega_t^{\nu}\rangle + \\
&+\langle \partial^k \mathsterling_i\omega_t^{\nu}, \partial^k
\mathsterling_i\omega_t^{\nu} \rangle \\
&=\langle\partial^k\mathsterling_i\omega_t^{\nu},L_i\omega_t^{\nu}\rangle +
\langle(L_i\mathsterling_i)\omega_t^{\nu}, \partial^k\omega_t^{\nu}\rangle. \end{aligned}
\end{equation*}

\noindent Due to our initial assumptions on the vector fields $(\xi_i)_i$, each term can be bounded by $%
\|\omega_t^{\nu}\|_{k,2}^2$, so the required inequality is proven. \newline\\
v. For $\beta \geq 2$ we have 
\begin{equation*}
\begin{aligned} &|\langle \partial^k \omega_t^{\nu_n, R, n}, K_R(\omega_t^{\nu_{n-1}, R, n-1})\partial^k(\mathsterling_{u_t^{\nu_{n-1},R,n-1}}
\omega_t^{\nu_n,R,n}) \rangle| \\
&=\big|\langle \partial^k 
\omega_t^{\nu_n,R,n}, K_R(\omega_t^{\nu_{n-1}, R,n-1})\displaystyle\sum_{|\beta|\leq k}
C_k^{\beta}(\partial^{\beta}u_t^{\nu_{n-1},R,n-1}) \cdot (\partial^{k-\beta}(\nabla
\omega_t^{\nu_n,R,n})) \rangle\big| \\ & =\big|K_R(\omega_t^{\nu_{n-1}, R, n-1})\displaystyle\sum_{|\beta|\leq k}
C_k^{\beta} \langle \partial^k \omega_t^{\nu_n,R,n}, \partial^{\beta}u_t^{\nu_{n-1},R,n-1}
\cdot \partial^{k-\beta}(\nabla \omega_t^{\nu_n,R,n}) \rangle \big|\\ & =
\bigg|K_R(\omega_t^{\nu_{n-1}, R, n-1})\displaystyle\sum_{|\beta|\leq k} C_k^{\beta}
\displaystyle\int_{\mathbb{T}^2} \partial^k \omega_t^{\nu_n,R,n}
\cdot\partial^{\beta}u_t^{\nu_{n-1},R,n-1} \cdot \partial^{k-\beta}(\nabla
\omega_t^{\nu_n,R,n})dx \bigg| \\ & \leq K_R(\omega_t^{\nu_{n-1}, R, n-1})\displaystyle\sum_{|\beta|\leq k}
C_k^{\beta} \displaystyle\int_{\mathbb{T}^2} |\partial^k \omega_t^{\nu_n,R,n}
\cdot\partial^{\beta}u_t^{\nu_{n-1},R,n-1} \cdot \partial^{k-\beta}(\nabla
\omega_t^{\nu_n,R,n})|dx \end{aligned}.
\end{equation*}
Using H\"older and Cauchy-Schwartz inequalities one has 
\begin{equation*}
\begin{aligned} \displaystyle\int_{\mathbb{T}^2} |\partial^k \omega_t^{\nu_n,R,n}
\cdot\partial^{\beta}u_t^{\nu_{n-1},R,n-1} \cdot \partial^{k-\beta}(\nabla
\omega_t^{\nu_n,R,n})|dx & \leq
\|\partial^k\omega_t^{\nu_n,R,n}\|_2\|\partial^{\beta}u_t^{\nu_{n-1},R,n-1}\cdot
\partial^{k-\beta}(\nabla \omega_t^{\nu_n,R,n})\|_2 \\ & \leq
\|\partial^k\omega_t^{\nu_n,R,n}\|_2\|\partial^{\beta}u_t^{\nu_{n-1},R,n-1}\|_4\|\partial^{k-%
\beta}(\nabla \omega_t^{\nu_n,R,n})\|_4 \\ &\leq
\|\omega_t^{\nu_n,R,n}\|_{k,2}\|\partial^{\beta}u_t^{\nu_{n-1},R,n-1}\|_4\|\partial^{k-\beta}(%
\nabla \omega_t^{\nu_n,R,n})\|_4 \end{aligned}
\end{equation*}
By Gagliardo-Nirenberg inequality 
\begin{equation*}
\begin{aligned} \|\partial^{\beta}u_t^{\nu_{n-1},R,n-1}\|_4 &\leq
C\|\partial^{\beta+1}u_t^{\nu_{n-1},R,n-1}\|_2^a\|u_t^{\nu_{n-1},R,n-1}\|_2^{1-a} \\ & \leq
C\|u_t^{\nu_{n-1},R,n-1}\|_{\beta+1,2}^a\|u_t^{\nu_{n-1},R,n-1}\|_2^{1-a} \\ & \leq C
\|u_t^{\nu_{n-1},R,n-1}\|_{k+1,2}^a\|u_t^{\nu_{n-1},R,n-1}\|_2^{1-a} \\ & \leq C
\|\omega_t^{\nu_{n-1},R,n-1}\|_{k,2}^a\|u_t^{\nu_{n-1},R,n-1}\|_2^{1-a} \\ & \leq C
\|\omega_t^{\nu_{n-1},R,n-1}\|_{k,2} \\ \end{aligned}
\end{equation*}
since $\|u_t^{\nu_{n-1},R,n-1}\|_2^{1-a} \leq \|u_t^{\nu_{n-1},R,n-1}\|_{2,2}^{1-a} \leq
\|u_t^{\nu_{n-1},R,n-1}\|_{k+1,2}^{1-a} \leq C\|\omega_t^{\nu_{n-1},R,n-1}\|_{k,2}^{1-a}$ by the
Biot-Savart law.\newline
\newline
Similarly 
\begin{equation*}
\begin{aligned} \|\partial^{k-\beta}(\nabla \omega_t^{\nu_n,R,n})\|_4 & \leq
C\|\partial^{k-\beta+1}(\nabla \omega_t^{\nu_n,R,n})\|_2^a
\|\nabla\omega_t^{\nu_n,R,n}\|_2^{1-a} \\ & \leq
C\|\nabla\omega_t^{\nu_n,R,n}\|_{k-\beta+1,2}^a \|\omega_t^{\nu_n,R,n}\|_{2,2}^{1-a} \\
& \leq C\|\omega_t^{\nu_n,R,n}\|_{k,2}^a\|\omega_t^{\nu_n,R,n}\|_{k,2}^{1-a} \\ & =
C\|\omega_t^{\nu_n,R,n}\|_{k,2} \\ \end{aligned}
\end{equation*}
Therefore 
\begin{equation*}
\begin{aligned} 
|\langle \partial^k \omega_t^{\nu_n,R,n},K_R(\omega_t^{\nu_{n-1},R,n-1})\partial^k(\mathsterling_{u_t^{\nu_{n-1},R,n-1}}
\omega_t^{\nu_n,R,n}) \rangle| &\leq K_R(\omega_t^{\nu_{n-1},R,n-1})\|\omega_t^{\nu_{n-1},R,n-1}\|_{k,2}\|\omega_t^{\nu_n,R,n}\|_{k,2}^2 \\
& \leq \tilde{C}\|\omega_t^{\nu_n,R,n}\|_{k,2}^2.
\end{aligned} 
\end{equation*}

\noindent For $\beta = 0:$ 
\begin{equation*}
\begin{aligned} \displaystyle\int_{\mathbb{T}^2}\partial^k\omega_t^{\nu_n,R,n} \cdot u_t^{\nu_{n-1},R,n-1}\cdot
\partial^k(\nabla\omega_t^{\nu_n,R,n})dx & = \displaystyle\int_{\mathbb{T}^2}
 u_t^{\nu_{n-1},R,n-1}\cdot \nabla(\partial^k\omega_t^{\nu_n,R,n}) \cdot \partial^k\omega_t^{\nu_n,R,n}dx \\ 
 & = \frac{1}{2}\displaystyle\int_{\mathbb{T}^2} u_t^{\nu_{n-1},R,n-1} \cdot \nabla |\partial^k\omega_t^{\nu_n,R,n}|^2dx\\
 & = -\frac{1}{2}\displaystyle\int_{\mathbb{T}^2} (\partial^k\omega_t^{\nu_n,R,n})^2(\nabla \cdot u_t^{\nu_{n-1},R,n-1})dx \\ & = 0. \end{aligned}
\end{equation*}
For $\beta = 1:$ 
\begin{equation*}
\begin{aligned} \|\nabla u_t^{\nu_{n-1},R,n-1}\cdot \partial^k\omega_t^{\nu_n,R,n}\|_2 & \leq
\|\nabla u_t^{\nu_{n-1},R,n-1}\|_{\infty}^2\|\partial^k\omega_t^{\nu_n,R,n}\|_2 \\ & \leq
C\|\nabla u_t^{\nu_{n-1},R,n-1}\|_{k,2} \|\omega_t^{\nu_n,R,n}\|_{k,2} \\ & \leq
C\|u_t^{\nu_{n-1},R,n-1}\|_{k+1,2} \|\omega_t^{\nu_n,R,n}\|_{k,2} \\ & \leq
C\|\omega_t^{\nu_{n-1},R,n-1}\|_{k,2}\|\omega_t^{\nu_n,R,n}\|_{k,2}. \end{aligned}
\end{equation*}
We used the embedding $\mathcal{W}^{k,2} \hookrightarrow L^{\infty}$ and the
Biot-Savart law. 
We need this property for relative compactness (see Proposition \ref{tightnessI}). \\
vi. After applying the It\^{o} formula we obtain 
\begin{equation*}
\begin{aligned}
\|\partial^k\omega_t^{\nu_n, R, n}\|_2^2 &= \|\partial^k\omega_0^{\nu_n, R, n}\|_2^2 + 2\nu_n\displaystyle\int_0^t\langle \partial^k\omega_s^{\nu_n, R, n}, \partial^{k+2}\omega_s^{\nu_n, R, n}\rangle ds \\
&- 2\displaystyle\int_0^t \langle \partial^k\omega_s^{\nu_n, R, n}, K_R(\omega_s^{\nu_{n-1},R, n-1})\partial^k\mathsterling_{u_s^{\nu_{n-1}, R, n-1}}\omega_s^{\nu_n, R, n}\rangle ds \\
& + \displaystyle\sum_{i=1}^{\infty}\displaystyle\int_0^t \langle \partial^k\omega_s^{\nu_n,R,n}, \partial^k\mathsterling_i^2\omega_s^{\nu_n,R,n}\rangle ds \\ 
& + \displaystyle\sum_{i=1}^{\infty}\displaystyle\int_0^t\langle \partial^k\mathsterling_i\omega_s^{\nu_n, R, n}, \partial^k\mathsterling_i\omega_s^{\nu_n, R, n}\rangle ds \\
& - 2\displaystyle\sum_{i=1}^{\infty}\displaystyle\int_0^t\langle \partial^k\omega_s^{\nu_n, R, n}, \partial^k\mathsterling_i\omega_s^{\nu_n, R, n} \rangle dW_s^{i,n}. \\
\end{aligned}
\end{equation*}
We analyse each term. One can write
\begin{equation*}
\langle \partial^k\omega_s^{\nu_n,R,n},\partial^{k+2}\omega_s^{\nu_n,R,n}\rangle = -\|\partial^{k+1}\omega_s^{\nu_n,R,n}\|_2^2 \leq 0. 
\end{equation*}
We want to estimate the other terms independently of $\nu_n$. All terms are estimated above. 
Let 
\begin{equation*}
B_t := \displaystyle\sum_{i=1}^{\infty}\displaystyle\int_0^t\langle \partial^k\omega_s^{\nu_n, R, n}, \partial^k\mathsterling_i\omega_s^{\nu_n, R, n} \rangle dW_s^{i,n} \ \ \ \ \hbox{and} \ \ \ \
\beta_t := \|\omega_s^{\nu_n, R, n}\|_{k,2}^2.
\end{equation*}
Obviously $B_t$ is a local martingale. Using the Burkholder-Davis-Gundy inequality there exists a constant $\tilde{\alpha}_p$ such that
\begin{equation*}
\mathbb{E}\big[\displaystyle\sup_{s\in[0,t]}|B_s|^p\big] \leq \tilde{\alpha}_p \mathbb{E}\left[\langle B\rangle_t^{\frac{p}{2}}\right].
\end{equation*}
Then 
\[
\displaystyle\mathbb{E}\big[\langle B\rangle_t^{\frac{p}{2}}\big] \leq C_{p,T}^2\int_{0}^{t}\mathbb{E}[\displaystyle\sup_{r\in[0,s]} \beta_r^p]ds
\]
following the same calculations as those from step v. and taking into account assumption \eqref{xiassumpt}. We have, for $t\in [0,T]$,
\begin{eqnarray*}
\beta_t &\leq& \beta_0 -2B_t+ (C_2'+C_1) \displaystyle\int_0^t\beta_sds \\
\sup_{s\in [0,t]}\beta_s^p&\le & 
3^{p-1} \beta_0^p + 2\times 3^{p-1}  \sup_{s\in [0,t]}| B_s|^p+ 3^{p-1}T(C_2'+C_1)\displaystyle\int_0^t\sup_{r\in [0,s]}\beta_r^{p}ds.
\end{eqnarray*}
In conclusion 
\begin{equation*}
\mathbb{E}\big[\displaystyle\sup_{s\in[0,t]} \beta_t^p\big] \leq \tilde{C_1}(T)\beta_0^p + \tilde{C_2}(T)\displaystyle\int_0^t \mathbb{E}\big[\displaystyle\sup_{s\in[0,t]} \beta_s^p]ds.
\end{equation*}
Finally,  using Gronwall's inequality, we deduce that
\begin{equation*}
\mathbb{E}\big[\displaystyle\sup_{s\in[0,t]}\|\omega_s^{\nu_n,R,n}\|_{k,2}^p\big] \leq C_p(T)
\end{equation*}
with $C(T):= \tilde{C}_1(T)\beta_0^p\exp(\tilde{C}_2(T)T)$. \\

The following lemma is instrumental in showing that the limit of the approximating sequence satisfies the Euler equation in $\mathcal{W}^{k,2}(\mathbb{T}^2)$ although the relative compactness property holds in $D\big([0,T], L^{2}(\mathbb{T}^2)\big)$. It is also essential when proving a priori estimates for the truncated solution $\omega^R$ (see Lemma \ref{aprioritruncatedEuler}). 
  
\begin{lemma}\label{lim}$\left.\right.$\\
i. Assume that $(a_n)_n$
is a sequence of functions such that $\displaystyle\lim_{n\mapsto\infty} a_n=a$ in $L^2(\mathbb{T}^2)$ and $\displaystyle\sup_{n>1}
\left\Vert a_n\right\Vert_{s,2}<\infty$ for $s \geq 0$. Then $a \in \mathcal{W}^{s,2}(\mathbb{T}^2)$ and $\left\Vert a \right\Vert_{s,2}<\displaystyle\sup_{n>1}\left\Vert a_n\right\Vert_{s,2}$. Moreover, $\displaystyle\lim_{n\mapsto\infty}a_n = a$ in $\mathcal{W}^{s',2}(\mathbb{T}^2)$ for any $s'<s$.\\
ii.  
Assume that $a_n:\Omega\mapsto \mathcal{W}^{s,2}(\mathbb{T}^2)$ is a sequence of measurable maps such that,  $\displaystyle\lim_{n\mapsto\infty} a_n=a$ in $L^2(\mathbb{T}^2)$, $\mathbb{P}$-almost surely or $\displaystyle\lim_{n\mapsto\infty} a_n=a$ in distribution. Further assume that $\displaystyle\sup_{n>1}\mathbb E[\left\Vert a_n\right\Vert_{s,2}^{p}]<\infty$. Then, $\mathbb{P}$-almost surely, $a\in \mathcal{W}^{s,2}(\mathbb{T}^2) $ and $\mathbb E[\left\Vert a\right\Vert_{s,2}^{p}] \le \displaystyle\sup_{n>1}\mathbb E[\left\Vert a_n\right\Vert_{s,2}^{p}]$, for any $p>0$. \\
\end{lemma}

\noindent\textbf{\textit{Proof of Lemma \ref{lim}}}\\
i. Since $\displaystyle\lim_{n \rightarrow \infty} a_n = a$ in $L^2(\mathbb{T}^2)$ it follows that for arbitrary $\lambda\in \mathbb{Z}^2$ we have 
\[
\lim_{n\mapsto\infty} \widehat{a_n}\left(  \lambda\right)=\lim_{n\mapsto\infty} \int_{\mathbb{T}^{2}}e^{2\pi i\lambda\cdot
x} a_n\left(  x\right)  dx=\int_{\mathbb{T}^{2}}e^{2\pi i\lambda\cdot
x} a\left(  x\right)  dx=\widehat{a}(\lambda)
\]
Therefore by Fatou's lemma
\begin{eqnarray*}
\left\Vert a\right\Vert _{s,2}^{2}= \sum_{\lambda\in\mathbb{Z}^{2}%
}(1+\left\vert \lambda\right\vert ^{2s})\left\vert
\widehat{a}\left(  \lambda\right)  \right\vert ^{2}&=&  \sum_{\lambda\in\mathbb{Z}^{2}%
}(1+\left\vert \lambda\right\vert ^{2s})\liminf_{n\mapsto\infty}\left\vert
\widehat{a_n}\left(  \lambda\right)  \right\vert ^{2}\\
&\le& \liminf_{n\mapsto\infty}\sum_{\lambda\in\mathbb{Z}^{2}%
}(1+\left\vert \lambda\right\vert ^{2s})\left\vert
\widehat{a_{n}}\left(  \lambda\right)  \right\vert ^{2}=\liminf_{n\mapsto\infty}\left\Vert a_n\right\Vert _{s,2}^{2}\le \sup_{n\ge 1}\left\Vert a_n\right\Vert_{s,2}^2.
\end{eqnarray*} 
For the second part we can write 
\begin{equation*}
\|a_n-a\|_{s',2} = \displaystyle\sum_{\lambda \in \mathbb{Z}^2, |\lambda|\leq M} (1+|\lambda|^{2s'})|(\hat{a}_n - \hat{a})\lambda| + \displaystyle\sum_{\lambda \in \mathbb{Z}^2, |\lambda|\geq M} (1+|\lambda|^{2s'})|(\hat{a}_n - \hat{a})\lambda|.
\end{equation*}
Note that 
\begin{equation*}
\begin{aligned}
\displaystyle\sum_{\lambda \in \mathbb{Z}^2, |\lambda|\geq M} (1+|\lambda|^{2s'})|(\hat{a}_n - \hat{a})\lambda| & \leq \displaystyle\sum_{\lambda \in \mathbb{Z}^2, |\lambda|\geq M} \frac{3(1+|\lambda|^{2s})}{1+M^{2(s-s')}}|(\hat{a}_n - \hat{a})\lambda| \\
& \leq \frac{3}{1+M^{2(s-s')}}\big(\displaystyle\sup_{n\geq 1}\|a_n\|_{s,2}^2 + \|a\|_{s,2}^2\big)
\end{aligned}
\end{equation*}
where the first inequality is true due to the fact that $|\lambda|\geq M$. Now we can choose $M$ such that the last term is strictly smaller than $\frac{\epsilon}{2}$. Likewise, $n$ can be chosen such that 
\begin{equation*}
\displaystyle\sum_{\lambda \in \mathbb{Z}^2, |\lambda|\leq M} (1+|\lambda|^{2s'})|(\hat{a}_n - \hat{a})\lambda| < \frac{\epsilon}{2}
\end{equation*}
hence $a_n$ converges to $a$ in $\mathcal{W}^{s',2}$. \\
ii. From above it follows that 
\[
\left\Vert a\right\Vert _{s,2}^{p}
=(\left\Vert a\right\Vert _{s,2}^{2})^\frac{p}{2}
\le ( \liminf_{n\ge 1}\left\Vert a_n\right\Vert_{s,2}^2)^\frac{p}{2}= \liminf_{n\ge 1}\left\Vert a_n\right\Vert_{s,2}^p
\]
and therefore
\[
\mathbb E [\left\Vert a\right\Vert _{s,2}^{p}]\le \liminf_{n\mapsto\infty}\mathbb E [\left\Vert a_n\right\Vert _{s,2}^{p}]\le \sup_{n\ge 1}\mathbb E [\left\Vert a_n\right\Vert_{s,2}^p].
\]

\begin{lemma}\label{aprioritruncatedEuler}
Let $\omega_t^R$ be the solution of the truncated Euler equation \eqref{truncatedEulere}.
There exists a constant $\tilde{C}(T)$ independent of $R$ such that 
\begin{equation*}
\mathbb{E}\big[\displaystyle\sup_{t \in [0,T]}\|\omega_t^{R}\|_{k,2}^4\big] \leq \tilde{C}(T). 
\end{equation*} 
\end{lemma}

\noindent Proof. We deduce that there exists a constant $\hat{C}(T)$  such that 
\begin{equation}\label{supoutside}
\sup_{t \in [0,T]}\mathbb{E}\big[\displaystyle\|\omega_t^{R}\|_{k,2}^4\big] \leq \hat{C}(T).  
\end{equation} from Lemma \ref{lim} ii. with 
\begin{equation*}
a_n := a_{n,t}:=\displaystyle\omega_t^{\nu_n,R,n} \ \ \ \hbox{and} \ \ \ a:=a_t:=\displaystyle\omega_t^R.
\end{equation*}
since $\displaystyle\lim_{n \rightarrow \infty}\omega_t^{\nu_n,R,n} = \omega_t^R$ $\mathbb{P}$ - almost surely in $L^2(\mathbb{T}^2)$ (see Section \ref{existenceEulertruncated}),
and 
$\displaystyle\sup_{n>1}\mathbb{E}[\|a_n\|_{s,2}^4] < \infty$
by Lemma \ref{apriori} vi. and Proposition \ref{tightnessI}. 

\begin{remark}\label{truncatedsimilar}
i. In Lemma \ref{aprioritruncatedEuler} we cannot follow an approach identical to the one used in Lemma \ref{apriori} vi, due to a lack of smoothness in the truncated equation \eqref{truncatedEulere}. This difficulty could be overcome if we embed $L^{2}\left(  \mathbb{T}^{2}\right)$ into  $L^{2}\left(  \mathbb{T}^{2};\mathbb{C}\right)$, consider the basis of functions $\left\{
e^{2\pi i\xi\cdot x};\xi\in\mathbb{Z}^{2}\right\} $, and express 
$ \omega_t^R$ as 
$
\sum_{\xi\in\mathbb{Z}^{2}%
}\left\vert \xi\right\vert ^{2k}\left\langle \omega_t^R,\varphi_\xi\right\rangle\varphi_\xi.
$
where $\varphi_\xi :\mathbb{T}^{2}\mapsto \mathbb C$, $\varphi_\xi(x)=e^{2\pi i\xi\cdot x}$, $x\in \mathbb{T}$ and $\xi\in\mathbb{Z}^{2}$.
Then $ \partial^{k}\omega_t^R$ will be square integrable if and only if $
\sum_{\xi\in\mathbb{Z}^{2}%
}\left\vert \xi\right\vert ^{2k}\left\langle \omega_t^R,\varphi_\xi
 \right\rangle ^{2}<\infty
$ and we can finish the proof using the weak form of the truncated equation. \\
ii. Except the estimates which involve the second order operator $\mathsterling_i^2$, all the other estimates derived in Lemma \ref{apriori} for $\omega^{\nu_n,R,n}$ hold also for $\omega^R$. 
\end{remark} 
In what follows we recall some basic results which have been used before. 

\begin{theorem}[Kurtz's criterion for relative compactness - 
\cite{EthierKurtz} Theorem 8.6]\label{kurtz} Let $(E, d)$ be a complete and
separable metric space, $(X^{\alpha})_{\alpha}$ a family of processes with càdlàg sample paths, and suppose that for every $\eta>0$ and any
rational $t\geq 0$ there exists a compact set $K_{\eta, t} \subset E$ such
that $\displaystyle\sup_{\alpha}\mathbb{P}\big(X_{t}^{\alpha} \notin
K_{\eta, t}\big) \leq \eta.$ Then the following two statements are
equivalent: \newline
a) $(X^{\alpha})_{\alpha}$ is relatively compact. \newline
b) For each $T>0$ there exists $\beta>0$ and a family $(\gamma_{\delta}^{%
\alpha})_{0<\delta<1, \ \hbox{all} \ \alpha}$ of nonnegative random
variables such that $\mathbb{E} \big[\tilde{d}^\beta (X_{t}^{\alpha},
X_{t+u}^{\alpha}) | \mathcal{F}_t^{\alpha}\big] \leq \mathbb{E}\big[%
\gamma_{\delta}^{\alpha} | \mathcal{F}_t^{\alpha} \big]
        $ and $\displaystyle \lim_{\delta \rightarrow 0}
\sup_{\alpha} \mathbb{E}\big[ \gamma_{\delta}^{\alpha} \big] =0 $ for $t \in
[0, T], u \in [0, \delta]$, where $\tilde{d}= d \wedge 1$ and the filtration $(\mathcal{F}_t^{\alpha})_t$ refers to the natural filtration $(\mathcal{F}_t^{X^{\alpha}})_t$.
\end{theorem}

\begin{theorem}{(Gagliardo-Nirenberg \cite{GagliardoNirenberg})}\label{gagliardonirenberg} \textit{Let $u \in L^{q}$ and $D^{m}u \in L^r$ with $1\leq q, r\leq
\infty$. Then there exists a constant C such that the following inequalities
hold for $D^{j}u$ with $0\leq j<m$:} 
\begin{equation*}
\|D^j u\|_p \leq C \|D^mu\|_r^a\|u\|_q^{1-a}
\end{equation*}
where $a \in [j/m, 1]$ is defined such that $\frac{1}{p} = \frac{j}{2} + a %
\big(\frac{1}{r}-\frac{m}{2}\big) + \frac{1-a}{q}.$
\end{theorem}


\begin{remark}\label{semigroup}
We denote by $(S^n(t))_t$ the semigroup of the generator $A:=\nu_n\Delta$.
This semigroup is strongly continuous (see \cite{Lunardi}) and for any $f \in L^2(\mathbb{T}^2)$
it is true that 
\begin{equation*}
\|S^n(t)f\|_{k,2} \leq \|f\|_{k,2}.
\end{equation*}
\end{remark}
\begin{theorem}{(Theorem 4.2 from \cite{KurtzProtterII})}\label{Kurtzstochint}
Let $(\mathcal{F}_t^n)_t$ be a filtration and $(X^n)_n$ a sequence of $(\mathcal{F}_t^n)_{t}$-adapted processes with c\`{a}dl\`{a}g trajectories. Let $(W^n)_n$ be a sequence of standard Brownian motions. If $(X^n, W^n)_n$ converges in distribution to $(X,W)$, in the Skorokhod topology, then $(X^n, W^n, \int X^n dW^n)$ converges in distribution to $(X, W, \int XdW)$ in the Skorokhod topology. If the first convergence holds in probability, then the convergence of the stochastic integrals holds also in probability.    
\end{theorem}

\end{document}